\documentclass[final,3p]{elsarticle}

\usepackage{lineno,hyperref}

\usepackage{amsmath,amsfonts,amssymb}
\usepackage{xfrac}
\usepackage{enumitem}
\usepackage{color}
\usepackage{graphicx}
\usepackage{caption}
\usepackage{subcaption}
\biboptions{sort&compress}

\begin{document}
	\makeatletter
	\def\ps@pprintTitle{%
		\let\@oddhead\@empty
		\let\@evenhead\@empty
		\let\@oddfoot\@empty
		\let\@evenfoot\@oddfoot
	}
	\makeatother

\begin{frontmatter}

\title{System Identification and  $H_\infty$-based Control of  Quadrotor Attitude}

\author[mainaddress]{Ali Noormohammadi-Asl \corref{correspondingauthor}}
\cortext[correspondingauthor]{Corresponding author}
\ead{anoormohammadi@ieee.org}

\author[mainaddress]{Omid Esrafilian }
\ead{esrafili@eurecom.fr}
\author[mainaddress]{Mojtaba Ahangar Arzati }
\ead{mojtaba.ahangar@email.kntu.ac.ir}
\author[mainaddress]{Hamid D. Taghirad}
\ead{Taghirad@kntu.ac.ir}
\address[mainaddress]{\textbf{A}dvanced \textbf{R}obotics and
	\textbf{A}utomated \textbf{S}ystems (ARAS)\fntext[arassite]{website: aras.kntu.ac.ir}\\Faculty of Electrical Engineering, K. N. Toosi University of Technology}

\begin{abstract}
The attitude control of a quadrotor is a fundamental problem, which has a pivotal role in a quadrotor stabilization and control. What makes this problem more challenging is the presence of uncertainty such as unmodelled dynamics and unknown parameters. In this paper, to cope with uncertainty, an  $H_\infty$ control approach is adopted for a real quadrotor. To achieve  $H_\infty$ controller, first a continuous-time system identification is performed on the experimental data to encapsulate a nominal model of the system as well as a multiplicative uncertainty. By this means,  $H_\infty$ controllers for both roll and pitch angles are synthesized. To verify the effectiveness of the proposed controllers, some real experiments and simulations are carried out. Results verify that the designed controller does retain robust stability, and provide a better tracking performance in comparison with a well-tuned PID and a $\mu$ synthesis controller. 
\end{abstract}


\end{frontmatter}


\section{Introduction}

Quadrotors are widely considered due to their remarkable features and applications in a variety of areas such as  monitoring \cite{wang2016detecting,berra2017commercial}, surveillance \cite{acevedo2014decentralized, kingston2008decentralized}, search and rescue \cite{silvagni2017multipurpose,khan2015information}, agriculture \cite {tokekar2016sensor,tetila2017identification}, and delivery of goods \cite{gawel2017aerial,arbanas2016aerial}.  
Quadrotors are also of paramount interest owing to the advantages of low cost, maneuverable, small size, simple design structure. 
The quadrotor is a nonlinear system with coupled states. In general, there are six states including the robot position $\left(x,y,z\right)$ and the angles (yaw, pitch, roll) with only four inputs, resulting in an underactuated system. Despite the mentioned advantages of quadrotors, the open-loop instability, coupled states, and underactuation make their control challenging. A quadrotor control is accompanied with several other serious challenges, particularly when it turns from a simple model in the theory to a practical problem in which many other key factors such as parameters uncertainty, unmodelled dynamics, and input constraints need to be considered. The main focus of this paper is the stabilization and attitude control of a quadrotor considering such uncertainties and input constraints.

Modeling and control of a quadrotor have been extensively studied in the literature. PID control is the most common controller widely used in commercial quadrotoros \cite{hoffmann2007quadrotor,szafranski2011different,garcia2012robust,zhou2017robust}. In \cite{szafranski2011different} three different PID structures are employed for the attitude control. In \cite{garcia2012robust} a robust PID controller is designed using affine parameterization. In \cite{zhou2017robust}  a cascade control for the attitude control is employed and a robust compensator is added to reduce the wind disturbance effect. Linear Quadratic Regulator (LQR) is another linear-based controller used for quadrotors \cite{bouabdallah2004pid,liu2013robust}. In LQR controllers, setting proper weights is difficult and needs trial and error. Linear-based controllers have simple designs and have been tested successfully on experimental and simulated platforms. Linear controllers, however, are based on a linear approximation of the actual system model by neglecting the nonlinear behavior of the system which may increase the probability of the system failure. Moreover, they are limited to low velocities as well as small angles of deviation, and their stability is guaranteed only near the selected operating point. 

To transcend the limitations of the linear controllers, many nonlinear control approaches have been proposed. 
In \cite{flores2013lyapunov}, a lyapunov-based controller is designed using singular perturbation theory for a quadrotor stabilization. Backstepping-based control  has also been developed by many researchers \cite{huo2014attitude,das2009backstepping,djamel2016attitude,shao2018robust}.   A backstepping controller using a lagrangian form of the quadrotor dynamics is designed in \cite{das2009backstepping}. In \cite{djamel2016attitude} an optimal backstepping controller using $\text{H}_\infty$  is proposed. Sliding mode control is another nonlinear control method which is capable of dealing with some classes of uncertainties and external disturbances \cite{chen2016robust,jia2017integral,yang2016attitude,zheng2014second,wang2019disturbance}. In \cite{chen2016robust,jia2017integral}, a combination of  the sliding mode control and backstepping control techniques are used for the control of quadrotors. Despite the mentioned capabilities, sliding mode control suffers from chattering, and it is necessary to know the upper bounds of uncertainties. 

In order to reduce the chattering effect, different approaches like fuzzy gain scheduling \cite{yang2016attitude} and higher order sliding mode have been proposed\cite{zheng2014second}. A fuzzy state observer is proposed in \cite{mallavalli2018fault} to estimate the unknown nonlinear functions of the uncertain system model, then an integral terminal sliding mode controller is used. Many other robust nonlinear control techniques have also  been proposed to cope with the quadrotor uncertainties. In \cite{liu2017robust}, a robust compensating technique is used to control quadrotors with time-varying uncertainties and delays. The design of   a nonlinear $\text{H}_\infty$ controller is studied in \cite{raffo2010integral} to achieve the robustness. In \cite{kerma2012nonlinear}, a nonlinear $\text{H}_\infty$ output feedback controller coupled with a high order sliding mode estimator is utilized to control the quadrotor in the presence of parameters uncertainty and external disturbances. The adaptive control is another  approach for controlling systems with unknown model parameters and have been widely used for controlling various practical systems\cite{hu2018adaptive,zhang2016active,jiang2019hydrothermal}. An adaptive controller is employed in \cite{tran2018adaptive} for trajectory tracking of a quadrotor with considering the input constraints and uncertain parameters in its nonlinear model.

In order to design a nonlinear-based controller, the quadrotor nonlinear model is required, which is hard to obtain. For instance, inertia matrix, center of gravity and some other parameters of a quadrotor are not easily accessible  and  it is difficult to obtain their accurate values. In most of the previous studies, a symmetric model for the quadrotor is assumed, however in real robots this assumption does not hold. Furthermore, some factors, like wind, cause environmental disturbances. Especially in an actual flight near obstacles or low attitude, undesirable wind effects on the stability and performance of a quadrotor is more obvious. Considering the wind model makes designing a nonlinear controller arduous and challenging. An approach to overcome the model uncertainty, such as unmodelled dynamics and parameter uncertainty, is using system identification techniques. In \cite{liu2018parameter}, a closed-loop multivariable extremum seeking algorithm (MESA) is suggested for parameter identification of a quadrotor. The black-box model identification using a continuous-time approach is adopted in \cite{bergamasco2014identification} to obtain a linear model for the dynamics of a quadrotor.

Linear robust control is a popular approach in the control theory that can help to mitigate the effect of unmodeled dynamics in linear-based controllers. This approach has been used in many practical problems such as surgery robots\cite{agand2017decentralized}, parallel robots \cite{bataleblu2016robust}, harmonic drive systems \cite{taghirad2001h}, and vibration control of elastic gantry crane  and 3D bar structure \cite{golovin2019robust,mystkowski2016mu}; and the results confirm its remarkable ability to deal with  uncertainty. In \cite{safaee2013system,wang2013robust}, the linear robust $\text{H}_\infty$ controller is employed for attitude and altitude control of a helicopter. In order to utilize the $\text{H}_\infty$ controller, a  comprehensive identification phase is required, which provides information about the linear model of the system and the associated uncertainty.

In this paper, there is no information available about the model and parameters (e.g. CAD model) of the quadrotor. Thus, in order to overcome the aforementioned limitations of linear and nonlinear controllers, embracing lack of system model and different sources of uncertainty, a linear robust controller has been derived and implemented on a real quadrotor. In this case, a continuous time system identification from the sampled data of the real robot is required to gather information about both the system model and uncertainty. The main contributions and phases of this paper are as follows:
\begin{itemize}
\item	The first phase is performing a system identification based on the experimental frequency response estimates to obtain a nominal linear model aligned with a multiplicative uncertainty block. This phase is the principle and most challenging part of control design because of system identification difficulties and its importance in designing an effective controller. In this phase, the first step is designing appropriate experiments which provide us with informative and usable data for system identification. Then, by applying existing tools and methods such as MATLAB identification and CONTSID toolbox, the nominal models and uncertainty weighting functions are acquired. 
\item Having obtained information about the system in the previous phase, a linear $H_\infty$ control is synthesized. In this phase, it is important to choose proper sensitivity and  input weighting  functions to achieve desired tracking and regulating performance.
\item Finally, the controller is implemented on the robot and its high performance and robustness are shown by simulation and experimental results. In addition, a PID and a $\mu$ analysis controllers are designed to compare their performance with $\text{H}_\infty$ controller.
\end{itemize}

\section{System and uncertainty encapsulation}\label{sec2}
In many systems, uncertainty emerges due to the dynamical perturbation, such as unmodelled and high frequency dynamics, in various parts of a system. To capture  uncertainty in the system, we utilize the classical multiplicative perturbation model. This model helps to encapsulate the various sources of uncertainty, e.g. parametric and unmodelled dynamics, in a full block of multiplicative uncertainty. Considering the nominal system transfer function as $G_0$, the family of uncertain systems may be shown as follows:
\begin{equation}\label{eq21}
\mathcal{G}=\left\{G\left(s\right) \mid G=\left(1+\Delta W\right)G_0 \right\} ,
\end{equation}
where $W$ is the uncertainty weighting function, and $\Delta$ is a stable perturbation satisfying $\lVert\Delta\rVert_\infty<1$.

Hence, for designing a robust control for a nonlinear system, an approximate linear model along with a weighting function providing information about the uncertainty profile, are required. 
In order to obtain a controller having proper performance in practice, it is necessary to perform a system identification on the experimental data generated from the real system. There are some challenges, however, for a quadrotor system identification:
\begin{itemize}
\item The quadrotor robots are inherent unstable,  thus a closed loop identification method is recommended. 
\item For the robust control design a continuous-time model is needed. 
\item The quadrotor system is a multi-input multi-output (MIMO) system, but due to its open-loop instability and piloting difficulty, the experimental data is better to be obtained for each channel separately, then perform the system identification. 
\end{itemize}
In what follows, we dwell on the closed loop system identification. For this, first we derive a theoretical model for the system, then identify the system based on the gathered experimental data. 

\subsection{Theoretical model}
\begin{figure} [h]
    \centering
    \includegraphics[height=2.5in]{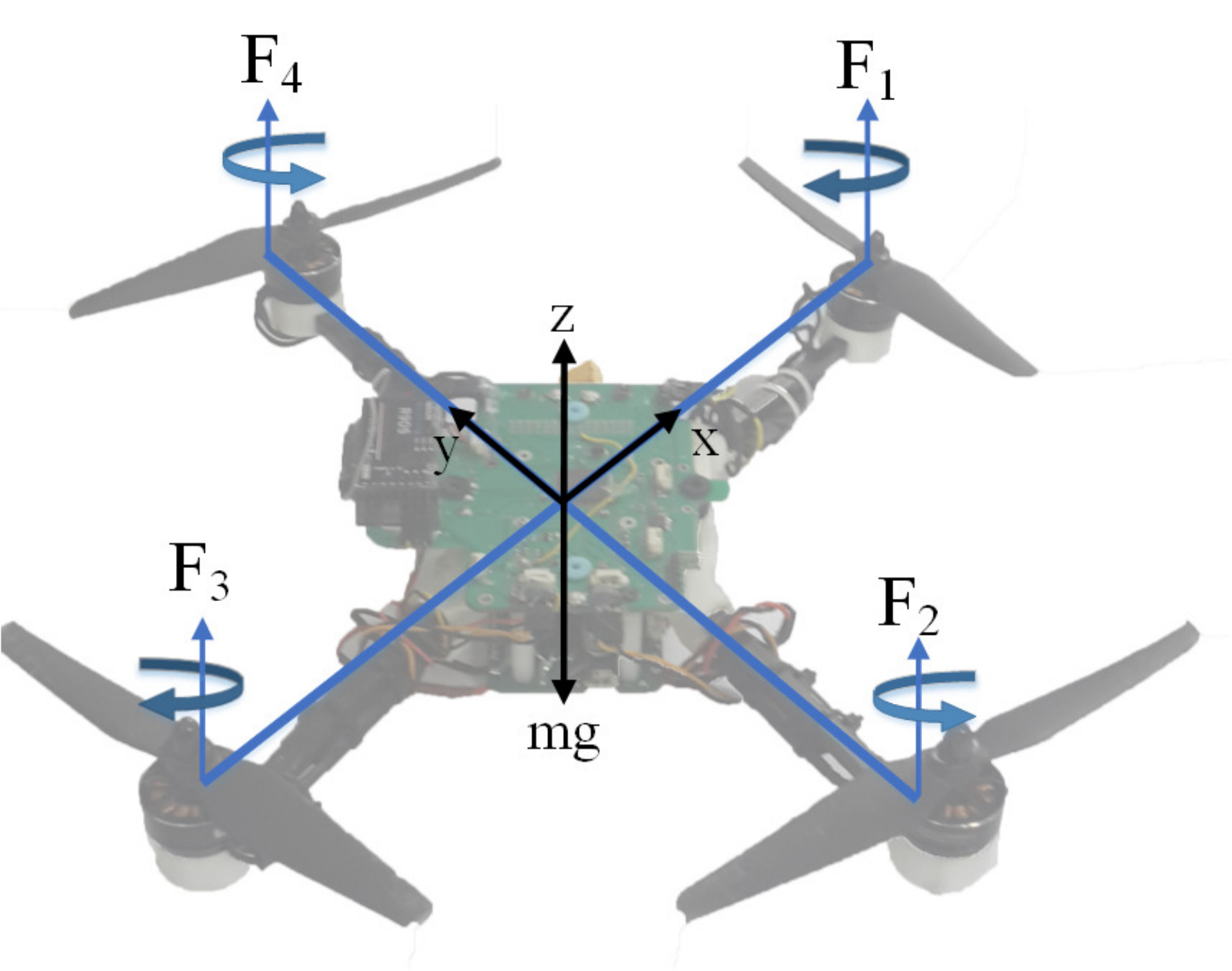}
    \caption{Quadrotor body frame, rotor arrangement and corresponding forces}
    \label{fig21}
\end{figure}

Modeling the dynamic of quadrotors has been studied extensively in the literature. Many models attempt to consider parameters and factors that are difficult to calculate or their values may vary in different situations such as the presence of wind. The goal of this paper is to design a controller that is largely independent of modeling parameters. 
To this end, first we model the dynamic of the robot. As shown in Fig. \ref{fig21} for a drone  with four rotors, which rotate in opposite directions (two rotors rotate in a clockwise direction and two other rotate counterclockwise), the angular accelerations pertain to the pitch, roll, and yaw angles are modeled by the following nonlinear models: 
\begin{equation}\label{eq22}
\begin{aligned}
&\ddot{\phi}=\frac{J_r\dot{\theta}\left(\Omega_1+\Omega_3-\Omega_2-\Omega_4\right)}{I_{xx}} +\frac{I_{yy}-I_{zz}}{I_{xx}}\dot{\psi}\dot{\theta}+\frac{bl\left(\Omega_2^2-\Omega_4^2\right)}{I_{xx}}, \\
&\ddot{\theta}=\frac{J_r\dot{\phi}\left(-\Omega_1-\Omega_3+\Omega_2+\Omega_4\right)}{I_{yy}} +\frac{I_{zz}-I_{xx}}{I_{yy}}\dot{\psi}\dot{\phi}+\frac{bl\left(\Omega_3^2-\Omega_1^2\right)}{I_{yy}}, \\
&\ddot{\psi}=\frac{d\left(\Omega_1^2+\Omega_3^2-\Omega_2^2-\Omega_4^2\right)}{I_{zz}} +\frac{I_{xx}-I_{yy}}{I_{zz}}\dot{\theta}\dot{\phi},
\end{aligned}
\end{equation}
where $\phi$, $\theta$, and $\psi$ stand for the roll, pitch and yaw angles, respectively. The remaining parameters are listed in table \ref{table21}. $\left(I_{yy}-I_{zz}\right)\dot{\psi}\dot{\theta}$, $\left(I_{zz}-I_{xx}\right)\dot{\psi}\dot{\phi}$ and $\left(I_{xx}-I_{yy}\right)\dot{\theta}\dot{\phi}$ are body gyro effects. Considering $\Omega_r=\left(\Omega_1+\Omega_3-\Omega_2-\Omega_4\right)$,  $J_r\dot{\theta}\Omega_r$ and $J_r\dot{\phi}\Omega_r$ represent propeller gyro effects. 

\begin{table}[t]
\caption{List of symbols in Eq. \ref{eq22}} \label{table21}
\begin{center}
\begin{tabular}{|c||c|}
  \hline
 Symbol  & Description \\
  \hline\hline
$I_{xx,yy,zz}$ & Inertia moments\\
\hline
$J_r$ & Rotor inertia\\
\hline
$\Omega_i$ & Propeller angular rate\\
\hline
$d$ & Drag factor\\
\hline
$b$ & Thrust factor\\
\hline
$l$ & Horizontal distance: propeller center to CoG{\textsuperscript{\tiny {*}}}\\
  \hline
\multicolumn{2}{l}{\small {Center of Gravity}}
\end{tabular}
\end{center}
\end{table}

\begin{figure} [h]
    \centering
    \includegraphics[height=2.0in]{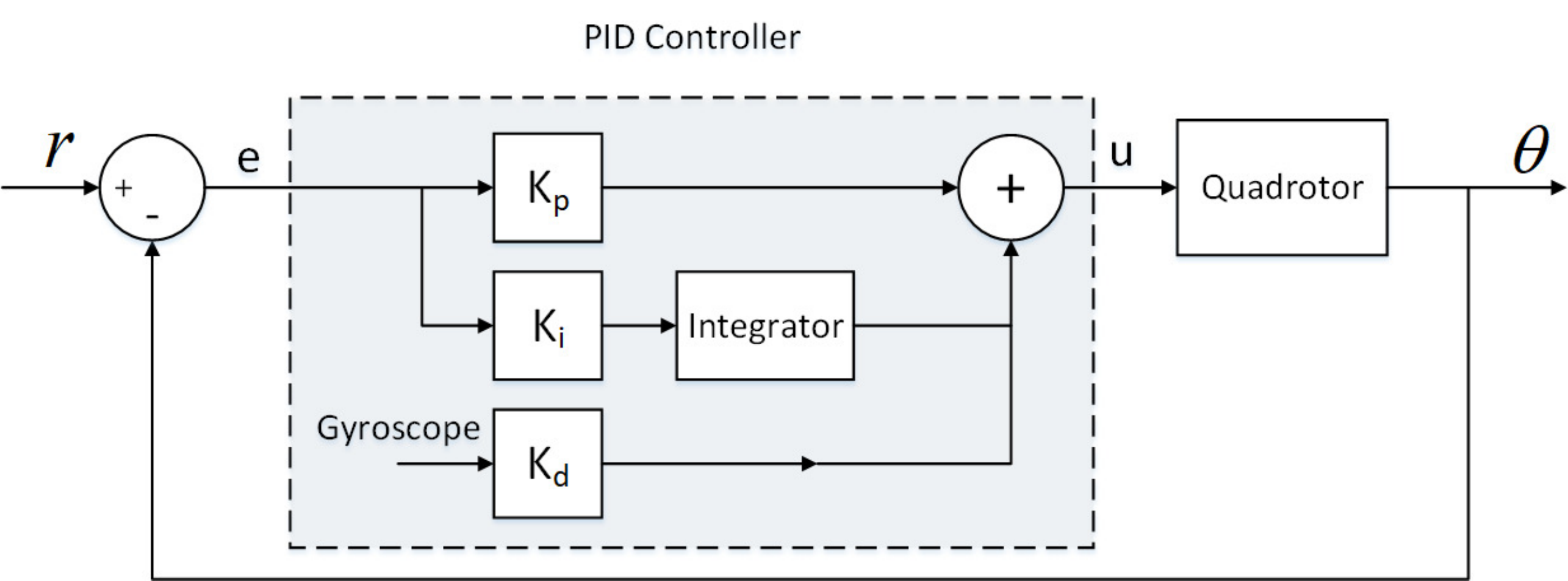}
    \caption{PID control structure for quadrotor}
    \label{fig22}
\end{figure}

In order to analyze the quadrotor model in the presence of a PID controller, the gyroscopic effects can be ignored compared to the motors' action. Thus, Eq. (\ref{eq22}) is rewritten as:
\begin{equation}\label{eq23}
\begin{aligned}
&\ddot{\phi}=\frac{bl\left(\Omega_2^2-\Omega_4^2\right)}{I_{xx}}, \\
&\ddot{\theta}=\frac{bl\left(\Omega_3^2-\Omega_1^2\right)}{I_{yy}}, \\
&\ddot{\psi}=\frac{d\left(\Omega_1^2+\Omega_3^2-\Omega_2^2-\Omega_4^2\right)}{I_{zz}}.
\end{aligned}
\end{equation}
The dynamics of the rotor is considered as $\frac{T_1}{s+T_2}$, in which $T_1,T_2\in \mathbb{R}^+$ and $T_2$ is the pole of the system and $\sfrac{T_1}{T_2}$ is the DC gain. Thus, the model of the robot is obtained as follows:
\begin{equation}\label{eq24}
\begin{aligned}
&\phi\left(s\right)=\frac{T_1^2bl}{s^2\left(s+T_2\right)^2I_{xx}}\left(u_4^2\left(s\right)-u_2^2\left(s\right) \right),\\
&\theta\left(s\right)=\frac{T_1^2bl}{s^2\left(s+T_2\right)^2I_{yy}}\left(u_3^2\left(s\right)-u_1^2\left(s\right) \right),\\
&\psi\left(s\right)=\frac{T_1^2d}{s^2\left(s+T_2\right)^2I_{zz}}\left(u_1^2\left(s\right)+u_3^2\left(s\right)-u_2^2\left(s\right) -u_4^2\left(s\right)\right), 
\end{aligned}
\end{equation}
where $u_1$, $u_2$, $u_3$, and $u_4$  are motors inputs.  Assuming the motors have fast response, the Eq. (\ref{eq24}) can be rewritten as follows, which is a double-integrator unstable system:
\begin{equation}\label{eq25}
\begin{aligned}
&\phi\left(s\right)=\frac{A_1}{s^2}U_1;\qquad U_1=u_4^2\left(s\right)-u_2^2\left(s\right),\\
&\theta\left(s\right)=\frac{A_2}{s^2}U_2; \qquad U_2=u_3^2\left(s\right)-u_1^2\left(s\right),\\
&\psi\left(s\right)=\frac{A_3}{s^2}U_3;\qquad U_3=u_1^2\left(s\right)+u_3^2\left(s\right)-u_2^2\left(s\right) -u_4^2\left(s\right),
\end{aligned}
\end{equation}
where $A_1, A_2$, and $A_3 \in \mathbb{R}^+$. The system model in the presence of a PID controller, as shown in Fig. \ref{fig22}, is obtained as:
\begin{equation}\label{eq26}
\begin{aligned}
&G_{roll}\left(s\right)=\frac{A_1\left(k_{d_r}s^2+k_{p_r}s+k_{i_r}\right)}{s^3+A_1\left(k_{d_r}s^2+k_{p_r}s+k_{i_r}\right)}U_1,\\
&G_{pitch}\left(s\right)=\frac{A_2\left(k_{d_p}s^2+k_{p_p}s+k_{i_p}\right)}{s^3+A_2\left(k_{d_p}s^2+k_{p_p}s+k_{i_p}\right)}U_2,
\end{aligned}
\end{equation}
where $K_{p_{r,p}}$, $K_{d_{r,p}}$, and $K_{i_{r,p}}$ are the nonnegative proportional, derivative, and integral gains, respectively. These equations give a useful insight into the number of poles and zeros of the system in the identification phase. 

\subsection{Control problem}
Attitude control is the integral part of the stabilization problem and helps to reach  the desired 3d orientation. In this paper, we focus on the roll and pitch angles  because they play the main role in the stability of the quadrotor. Due to the coupling of the system's states, the performance of roll and pitch has  considerable effect on other states. For example, when the quadrotor moves forward,  the pitch angle  changes and then returns to  zero. Thus, a fast and reliable control on roll and pitch angles is crucial. In addition, the disturbance exerts significant influence on the performance and stability of roll and pitch angles and consequently the overall performance of the robot. 

By considering $X=\left(\phi,\dot{\phi},\theta,\dot{\theta}\right)$, the control problem is to reach $X_d=\left(\phi_d,0,\theta_d,0\right)$ in which $\phi_d$ and $\theta_d$ are desired roll and pitch angles. In our case, the $\phi_d$ and $\theta_d$ are set to zero to stabilize the quadrotor. However, as mentioned before, overcoming model uncertainty is the main goal of this paper which should be achieved with considering two following conditions.
\begin{enumerate}
	\item	Having a high tracking performance and disturbance rejection
	
	To achieve this goal, we need that the controller provides a short settling time (about $0.3s$) and a small overshoot (about $10^{-6} rad$). 
	
	\item	Avoiding actuator saturation
	
	Due to the limitations of the motor, the control effort should be in a feasible range. In addition, actuator saturation may cause instability of the system. 
\end{enumerate}
Briefly, the problem is finding a controller providing the stability for the quadrotor by striking a balance  between robustness and performance of the system, in the presence of uncertainty and the actuator saturation.

\subsection{Experimental setup}
\begin{figure}
    \centering
    \begin{subfigure}[b]{0.48\textwidth}
        \includegraphics[height=0.8\linewidth,width=0.8\textwidth]{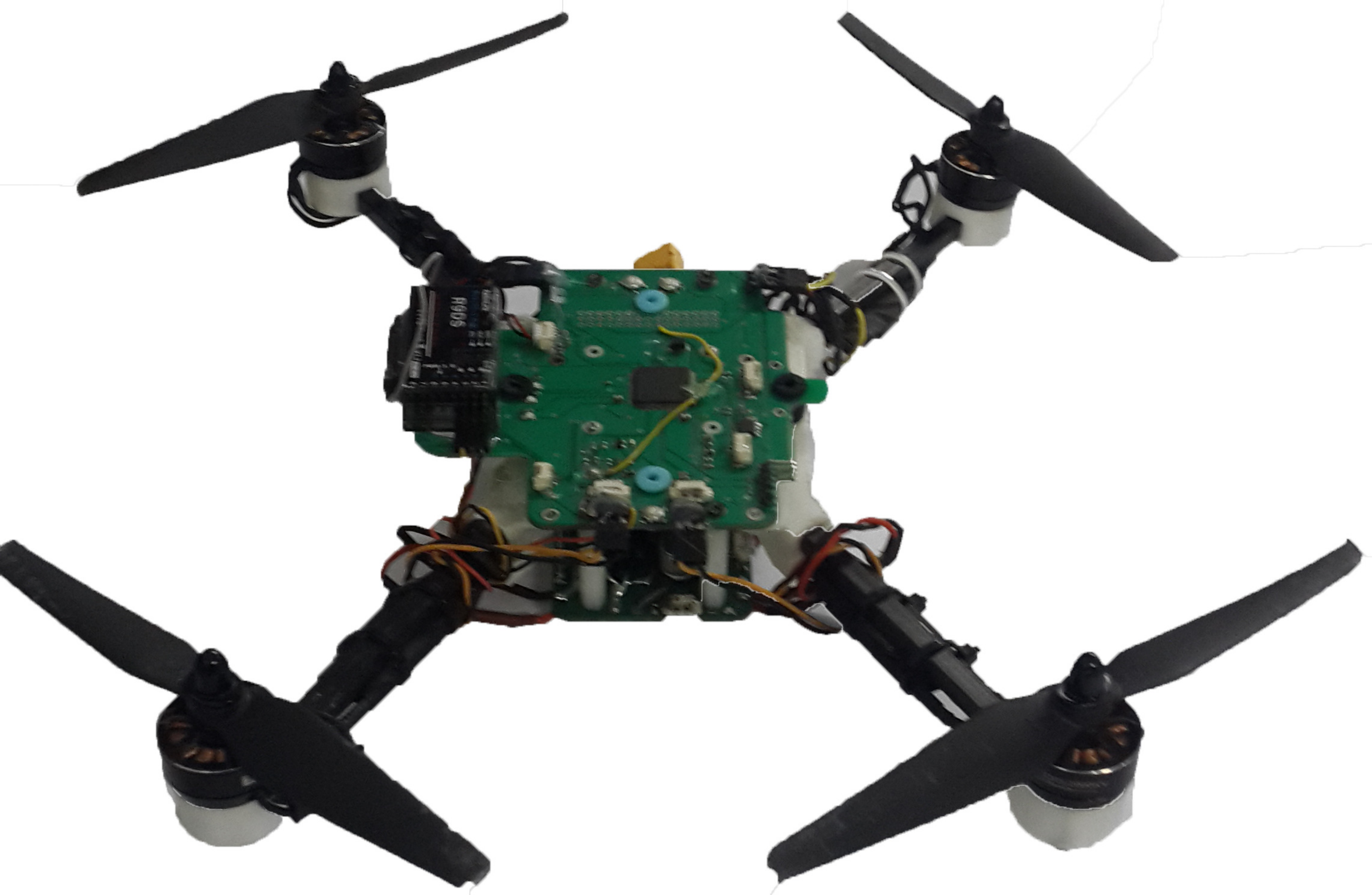}
        \caption{Quadrotor used in this paper}
        \label{fig2a}
    \end{subfigure}
    \begin{subfigure}[b]{0.48\textwidth}
        \includegraphics[height=0.8\linewidth,width=0.85\textwidth]{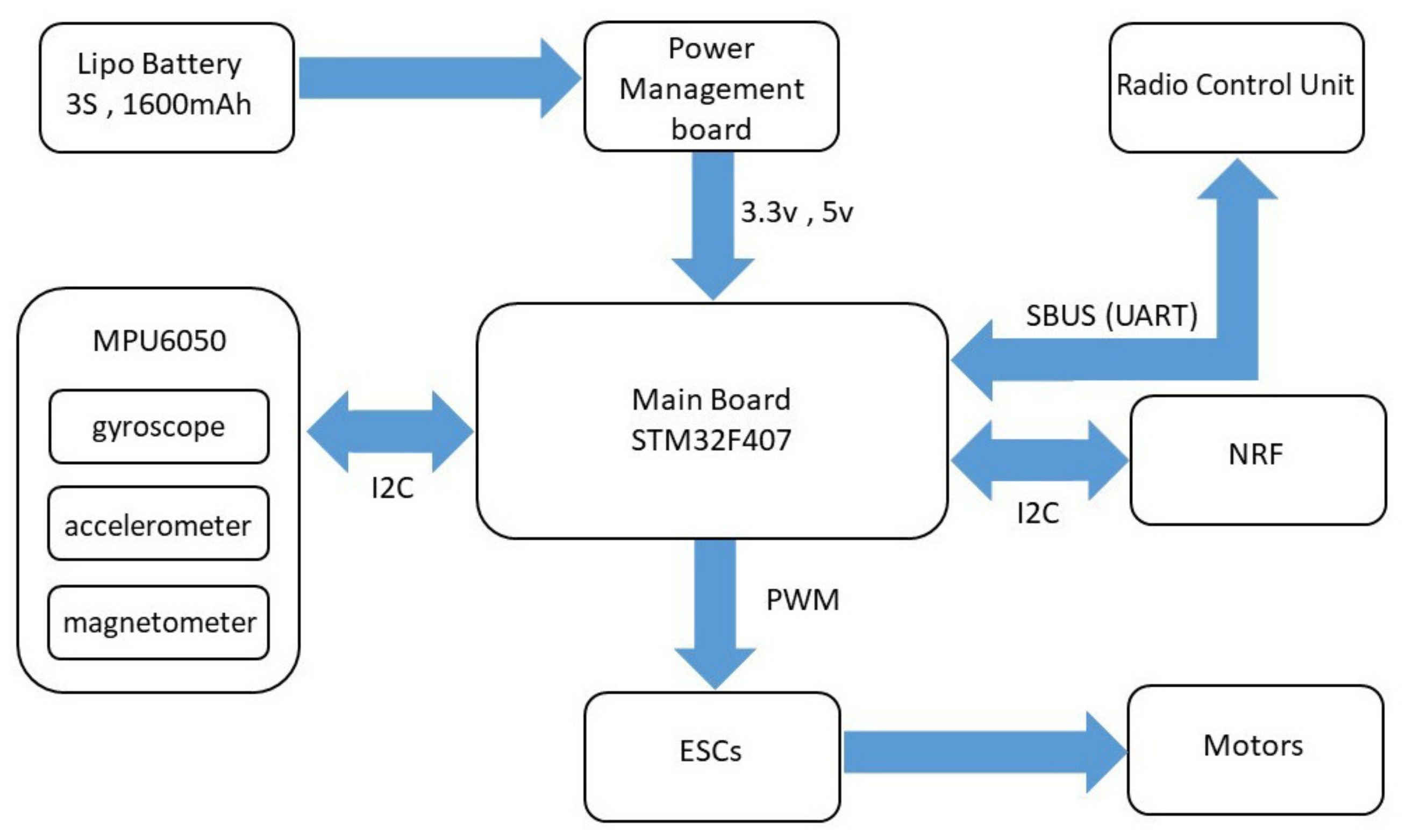}
        \caption{Block diagram of quadrotor’s units Interconnection}
        \label{fig2b}
    \end{subfigure}
    \caption{Experimental setup }\label{fig2}
\end{figure}

The quadrotor used in the experiments (Fig \ref{fig2a}) is based on a typical quadrotor design with some structural modifications. The power management circuit board and the carbon fiber tubes are used as the frame hub and motors arm, respectively. The main goal in the design of this quadrotor is achieving a compact integration of the mechanical structure and electronic devices in order to reduce the weight while maintaining a symmetric mass distribution with a center of gravity close to the center of the cross. The vehicle propeller-tip to propeller-tip distance is $38 \,cm$. This robot, embracing all the modules, sensors, battery, motors, and props has a mass of $570\,g$, and provides a total thrust of 850 g. The flight time is approximately 9 minutes with a 3-cell $1300\, mAh$ Li-Po battery. 


Three main units which play significant roles in the control of the robot are explained below, and the block diagram of their structure and interconnections are illustrated in Fig. \ref{fig2b}.

\textbf{\em{Embedded electronics:}}
The robot configuration consists of a flight control, Inertial Measurement Unit (IMU), radio control unit, telemetry module, and external IO. The flight control board, which is the main board of the robot, runs on an ARM STM32F407 micro-controller, which can operate up to $168\, MHz$. Other units interface with this controller via I2C and UART protocols. Each brushless motor has its own ESC (Electronic Speed Control) units, which is connected to the main board. PWM (Pulse-Width Modulation) is used to communicate with ESC at the speed of $400\, Hz$, and the main board sends control signals to ESC to adjust the rotational speed of the motor. 

\textbf{\em{Inertial Measurement Unit (IMU):}}
The IMU utilized in this quadrotor is MPU-6050, which contains a 3-axis gyroscope, 3-axis accelerometer, 3-axis magnetometer. The main board communicates with the MPU-6050 via I2C bus and runs a sensor fusion algorithm with the frequency of $500 \,Hz$. The quadrotor attitude is computed by using Mahony filter algorithm \cite{mahony2008nonlinear}, which estimates the absolute roll and pitch angles via fusing the acceleration vector and angular speeds measured by the gyroscopes. The angular speed which is measured by the gyroscope associated with yaw axis is used for estimating the yaw angle.

\textbf{\em{Wireless communication systems:}}
We established two types of system to support real time communications between the flying quadrotor and our personal computer as a ground control system (GCS): a digital radio telemetry unit ($915\,MHz$) and an analog radio link ($2.4\,GHz$).  We use NRF modules as a digital radio telemetry for data acquisition. The telemetry wireless interfaces with the main board using a UART serial protocol. Necessary flight data can be either saved on the internal memory of the micro-controller or sent to the GCS at rates of up to $115200\, bps$ using telemetry communication. Futaba remote controller is used as the analog radio link. The radio control (RC) receiver is connected to the main board through the SBUS protocol that allows to receive flight commands at the speed of $100 \,K bps$.

%

\subsection{System identification}
The problem of system identification of an aircraft or rotorcraft is a topic of interest and have been widely studied. In most of the studies, the system identification methods are used to obtain parameters value in the model of the system. However, in this paper, a black-box system identification, merely based on input and output experimental data, is performed. 

The first and basic step in the identification phase is selecting a proper input for applying to the system. In the frequency-domain identification, which is desired in this paper, a periodic excitation input is much preferred because it helps to cover and sweep the desired frequency spectrum effectively. The excitation inputs are better to be applied automatically (e.g programmed on the microcontroller) rather than manually (e.g radio controller). Thus, in order to obtain a continuous time model, a frequency-domain identification method is performed in which periodic chirp functions with different frequencies and amplitudes are implemented on the on-board micro-controller, and are used as exciting inputs. During the experiments, the identification for the pitch and roll angles are executed separately, and the quadrotor is supposed to track the reference chirp input. For example, if the chirp signal is selected as a reference input for the roll angle, the quadrotor uses the implemented simple PID controller for the roll angle to track the reference. Setting proper frequencies and amplitudes for chirp functions is important, because an improper frequency or amplitude, which are beyond the quadrotor ability to track, may cause inappropriate tracking and produces misleading information for identification. In this identification, the frequency and amplitude range vary in $\left[0.05 Hz,5 Hz\right]$ and $\left[5^\circ,11^\circ\right]$, respectively, which are feasible and reasonable pitch and roll angles.

\begin{figure}
    \centering
    \begin{subfigure}[b]{0.48\textwidth}
        \includegraphics[height=0.9\textwidth,width=\textwidth]{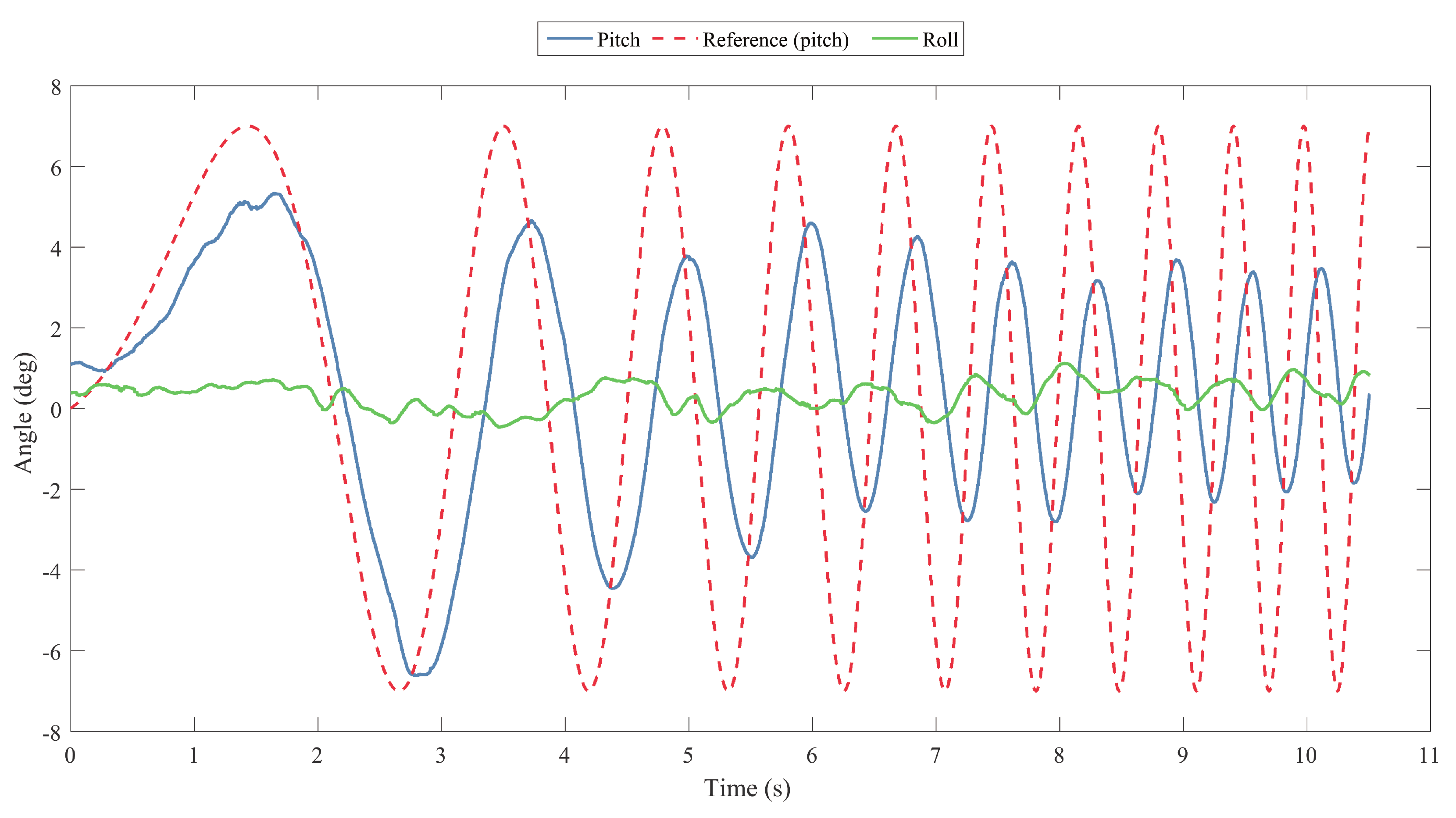}
        \caption{Coupling effect of the pitch on the roll}
        \label{fig231}
    \end{subfigure}
    ~ 
    \begin{subfigure}[b]{0.48\textwidth}
        \includegraphics[height=0.9\textwidth,width=\textwidth]{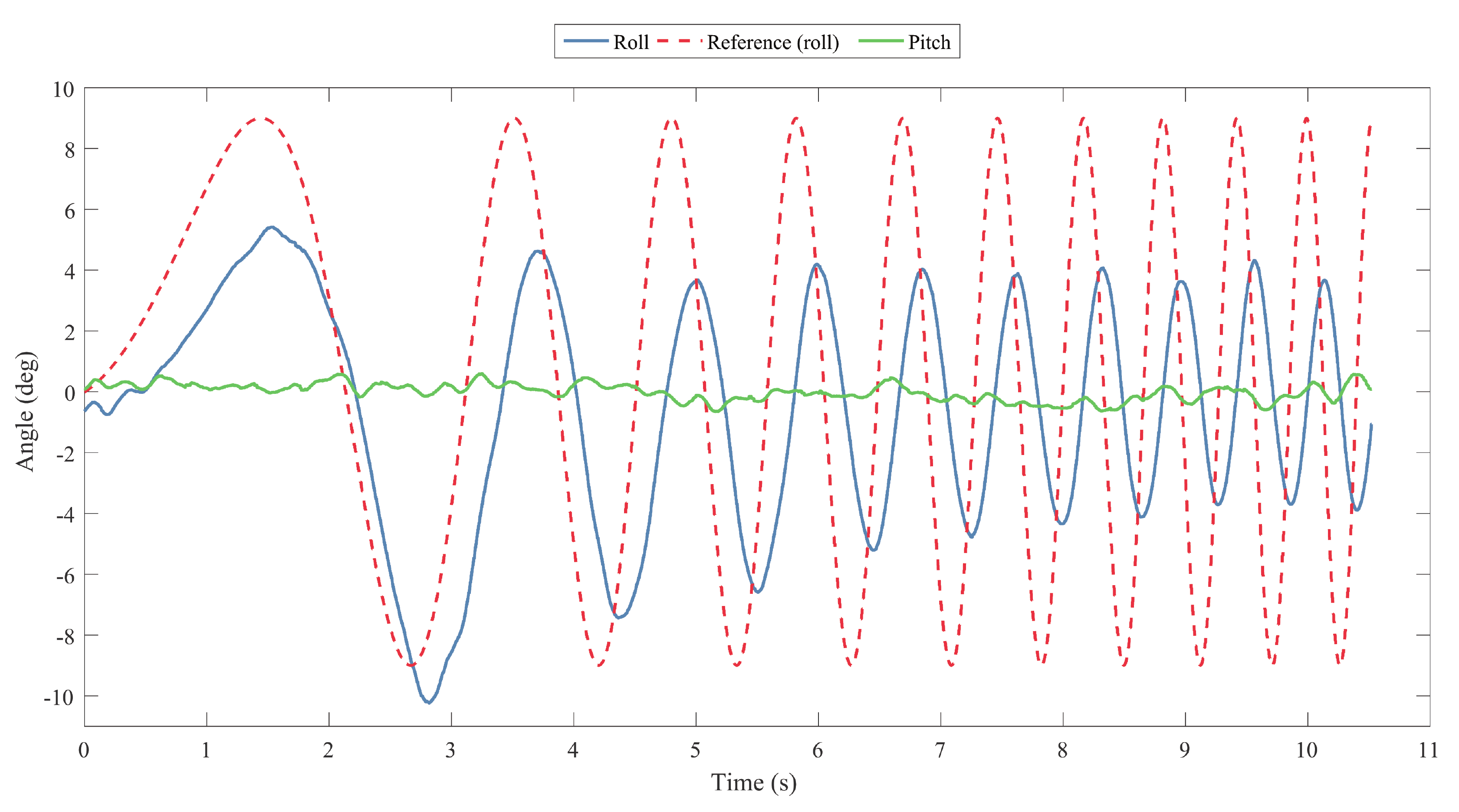}
        \caption{Coupling effect of the roll on the pitch}
        \label{fig232}
    \end{subfigure}
    \caption{The coupling effects}\label{fig23}
\end{figure}

Prior to the identification phase, the probable existence of delay is studied using correlation analysis method, and according to its results, the delay of the system is negligible. Secondly, the effect of exciting the roll angle on the pitch angle, and vice versa, is analyzed. In Fig. \ref{fig231}, the pitch angle, the desired reference pitch angle, and the roll angle are plotted, for one of the experiences. As shown in this figure, the coupling is negligible and the roll angle is almost independent of the pitch angle, thus their coupling can be considered as an unmodelled system uncertainty. The same condition holds for the pitch angle as well (Fig. \ref{fig232}). As a result, off-diagonal elements of the system model can be considered zero, which implies that the system can be represented by two SISO subsystems. Although this assumption may make the controller task more difficult at high frequencies, the controller design becomes much easier.

After performing several experiments for both pitch and roll angles, a linear transfer function should be fitted for each of them to obtain the stable closed-loop system model. We use MATLAB identification toolbox  for obtaining continuous-time transfer function using time-domain data.
In this step, it is important to select a proper number of poles and zeros. For this purpose, in each identification, we use different combinations of zeros and poles. However, to reduce the number of these combinations, the number of poles and zeros can be assumed close to  the approximate analytical transfer function of the system, obtained in Eq. (\ref{eq26}). Finally, for the model identification, the number of poles and zeros are considered in the ranges of $\left[2,5\right]$ and $\left[0,3\right]$, respectively. Consequently, for each experiment, fifteen  proper transfer functions are obtained, and one of them is selected as the best-fitted model. In order to select a transfer function, fitness percentage of different fitted models has a decisive role. However, there are two other factors which must be considered.
First, it is preferred that the selected transfer functions have  similar frequency responses, meaning that the zero/pole numbers should be close to each other. 
Second, the transfer functions with right hand-side zeros are not selected because the quadrotor system does not show the characteristics of a system with right half-plane zeros.


Based on the obtained transfer functions, for the pitch angle,  the real part of poles mainly lie in the range of $\left[-6,0\right]$, and  for the roll angle, real part of poles  are in the ranges of $\left[-7,0\right]$ and $\left[-10,-8\right]$. In both cases, the zeros lie in the range of $\left[-1,1\right]$. These zeros, which are located near the imaginary axis, mostly appear in pairs with poles very close to them, and therefore, they may be ignored due to the zero-pole cancellation.

%

\begin{figure}
    \centering
    \begin{subfigure}[b]{0.48\textwidth}
         \includegraphics[height=0.9\textwidth,width=\textwidth]{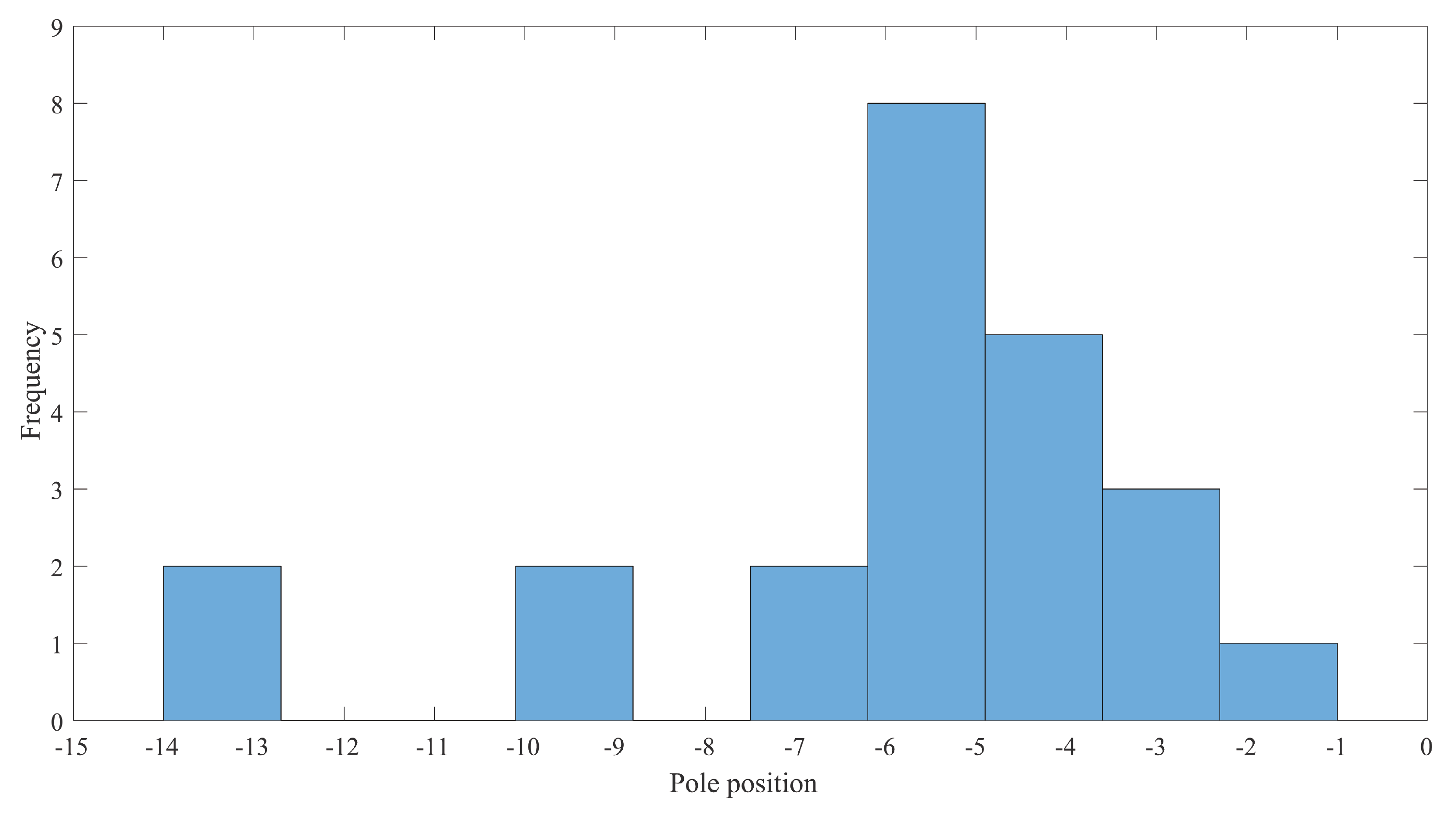}
        \caption{Pitch}
        \label{fig261}
    \end{subfigure}
    ~ 
    \begin{subfigure}[b]{0.48\textwidth}
         \includegraphics[height=0.9\textwidth,width=\textwidth]{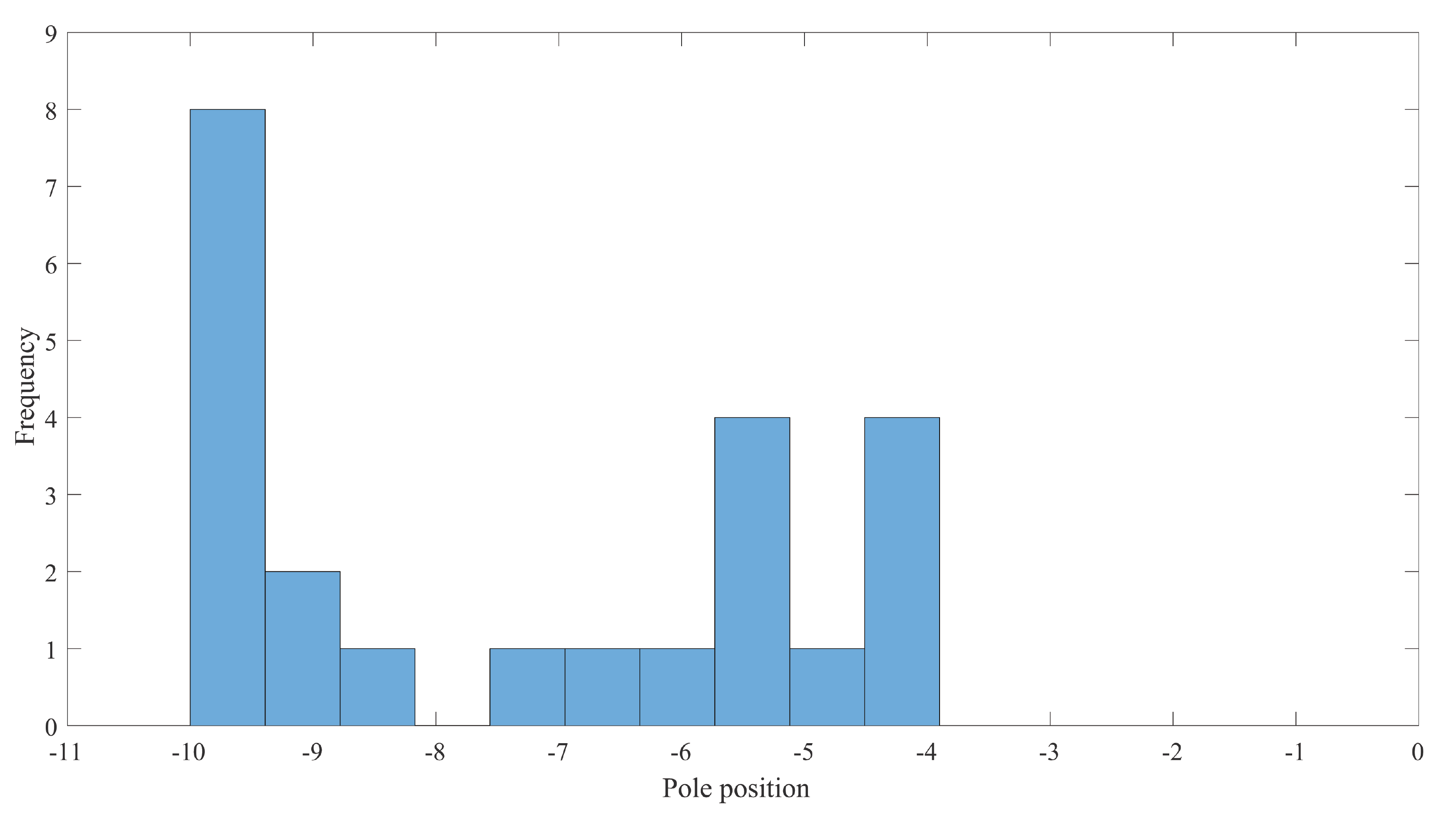}
        \caption{Roll}
        \label{fig262}
    \end{subfigure}
    \caption{Distribution of poles in transfer functions with 3 poles and no zero}\label{fig26}
\end{figure}

\begin{figure}
    \centering
    \begin{subfigure}[b]{0.48\textwidth}
         \includegraphics[height=0.9\textwidth,width=\textwidth]{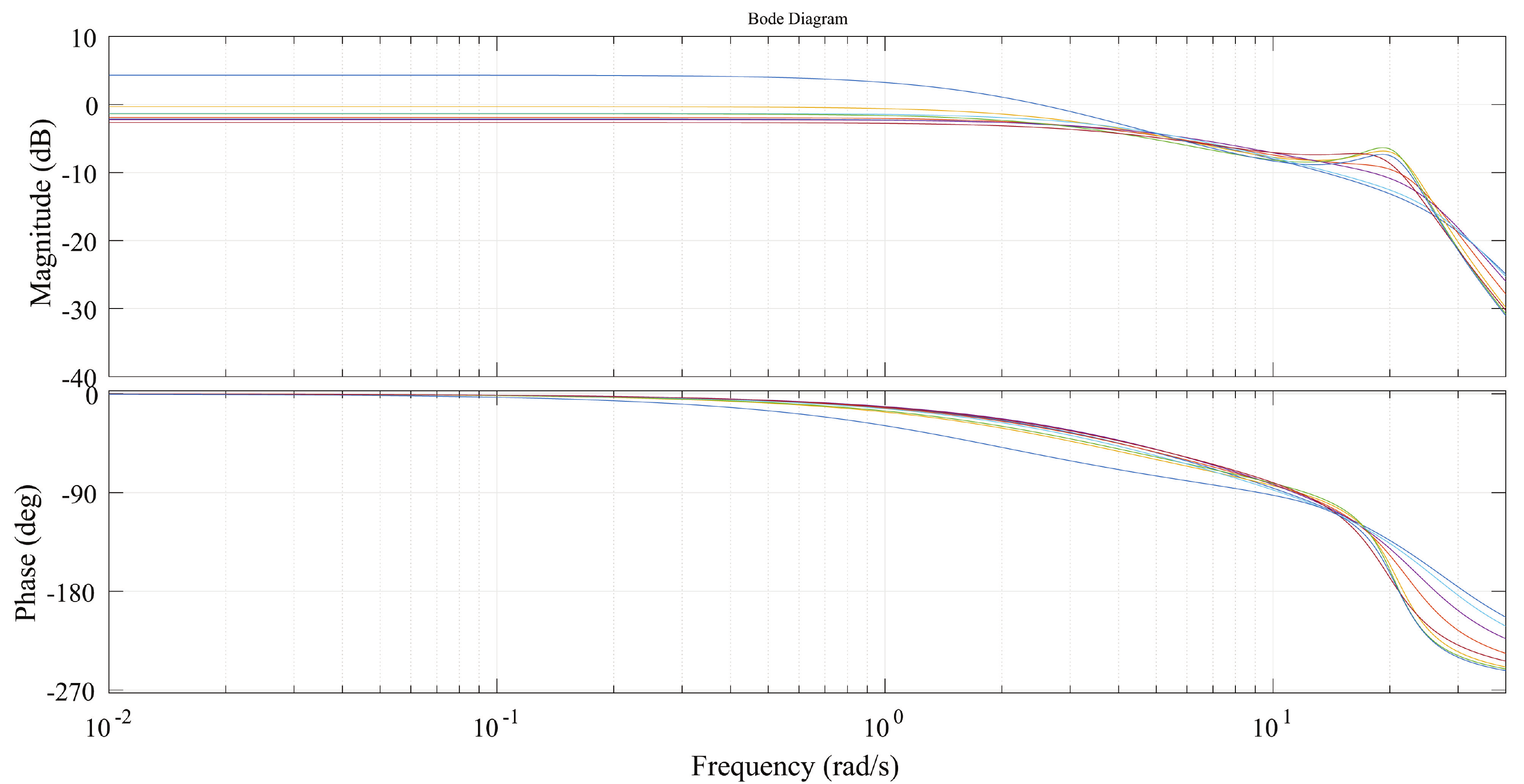}
        \caption{Pitch}
        \label{fig271}
    \end{subfigure}
    ~ 
    \begin{subfigure}[b]{0.48\textwidth}
         \includegraphics[height=0.9\textwidth,width=\textwidth]{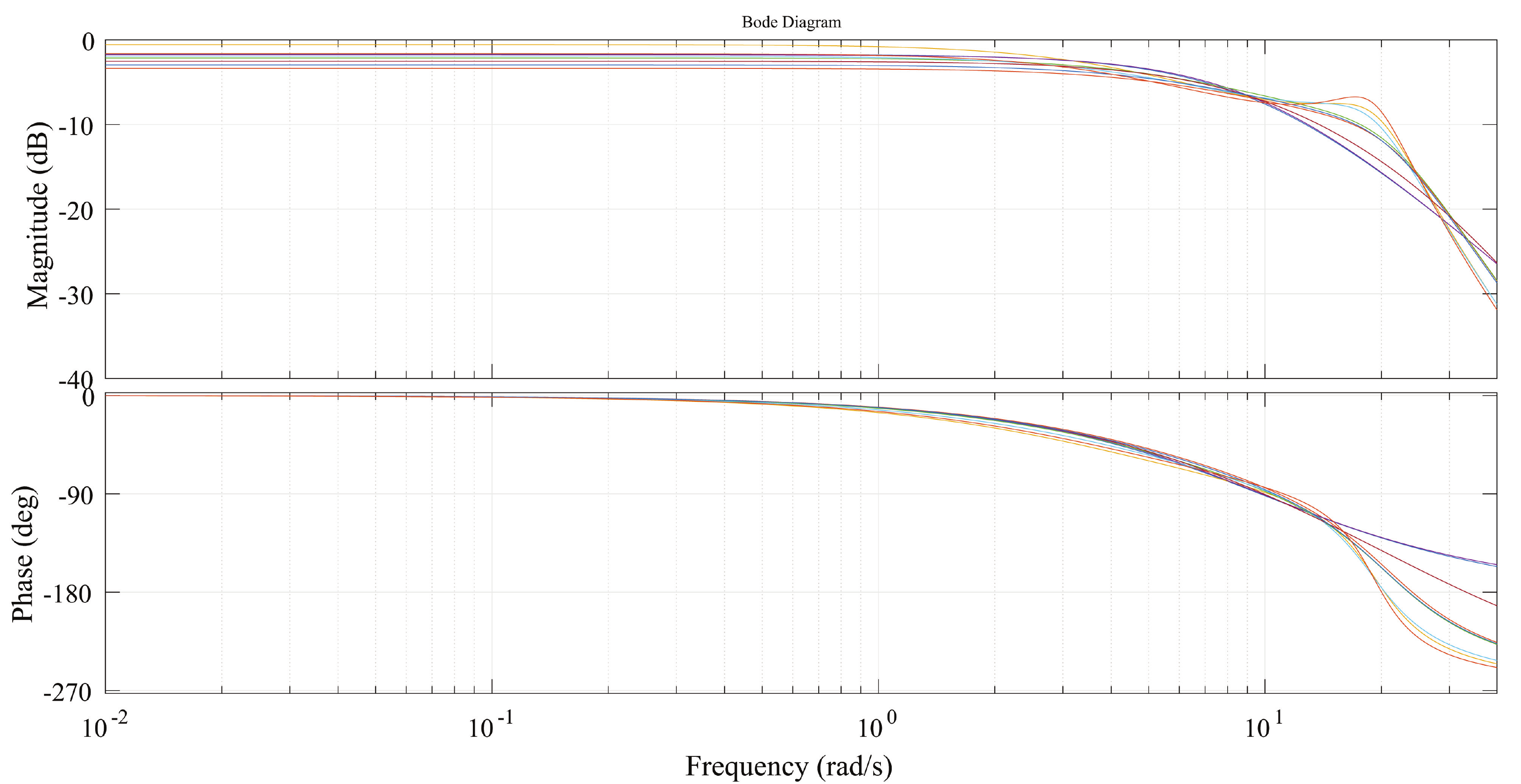}
        \caption{Roll}
        \label{fig272}
    \end{subfigure}
    \caption{Bode diagram of selected transfer functions in the experiments}\label{fig27}
\end{figure}

By analyzing the obtained transfer functions according to the three aforementioned criteria, it is observed that in the most cases, transfer functions with three poles and no zero can suitably describe the quadrotor system behavior. Figure \ref{fig261} and \ref{fig262} show the histogram of the real part of poles in the transfer functions of the pitch and roll angles, with three poles and no zero, respectively. According to these figures, it is expected that poles lie in the range of $\left[-7, -3\right]$ and $\left[-10, -8\right]$. 
In Fig. \ref{fig27}, the Bode diagram of the transfer functions obtained for the pitch and roll angle are depicted. It is clear from Fig. \ref{fig27} that the behavior of the transfer functions is similar, especially in frequencies less than $15 Hz$ where the robot mostly operates in. Now, by having the system transfer functions pertain to both axes of pitch and roll, we select one of them as the nominal transfer function. Hence, the transfer function which  lies near to the median of all responses, is considered as the nominal model.
The selected nominal model transfer function for the pitch and roll angle are as follows:
\begin{align}
G_{pitch}\left(s\right)=\frac{1547.4}{\left(s+5.373\right)\left(s^2+10.12s+390.4\right)},\label{eq27}\\
G_{roll}\left(s\right)=\frac{2049.8}{\left(s+6.764\right)\left(s^2+19.03s+426.2\right)}.
\end{align}

\begin{figure}
    \centering
    \begin{subfigure}[b]{0.48\textwidth}
         \includegraphics[height=1.1\textwidth,width=\textwidth]{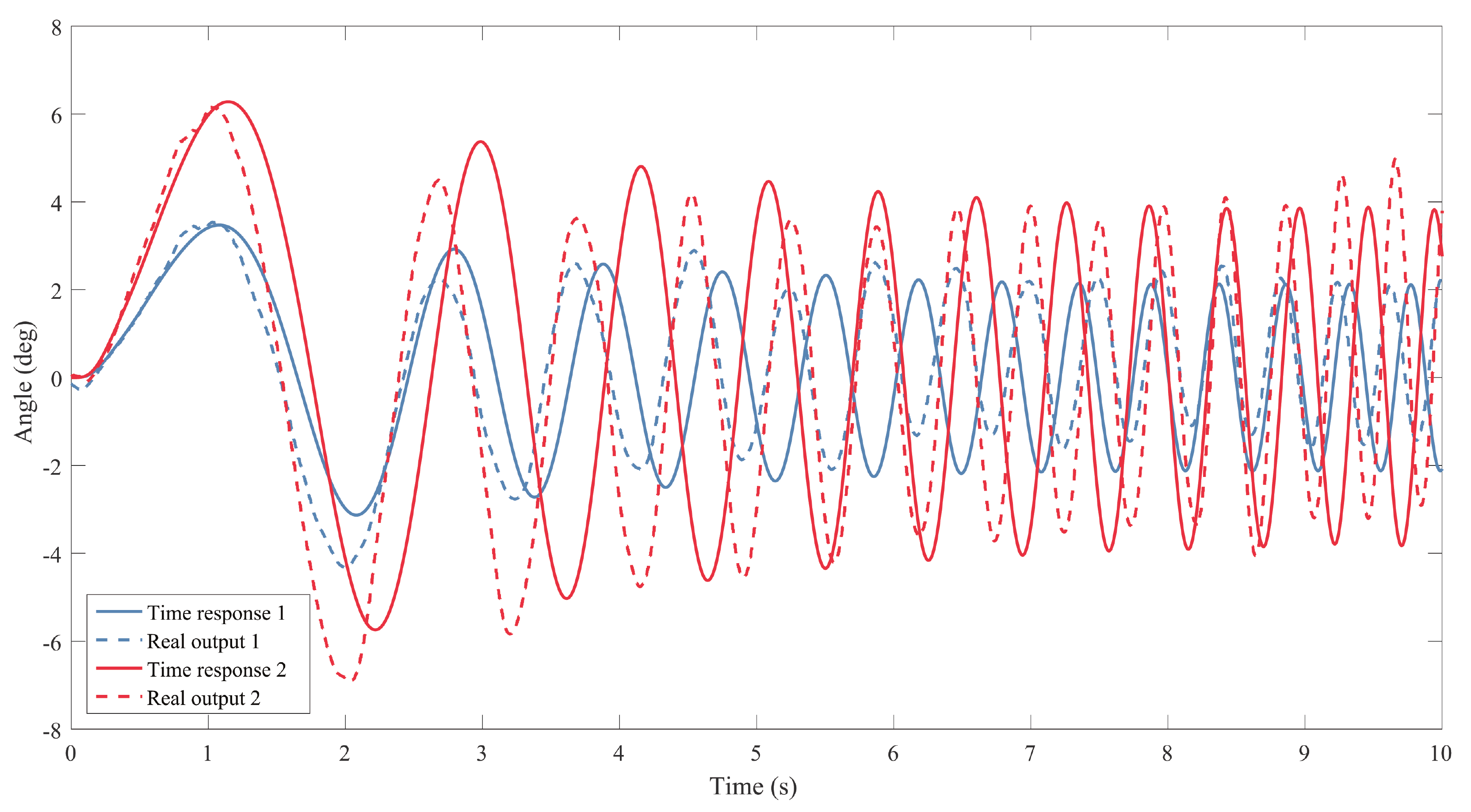}
        \caption{Pitch}
        \label{fig281}
    \end{subfigure}
    ~ 
    \begin{subfigure}[b]{0.48\textwidth}
         \includegraphics[height=1.1\textwidth,width=\textwidth]{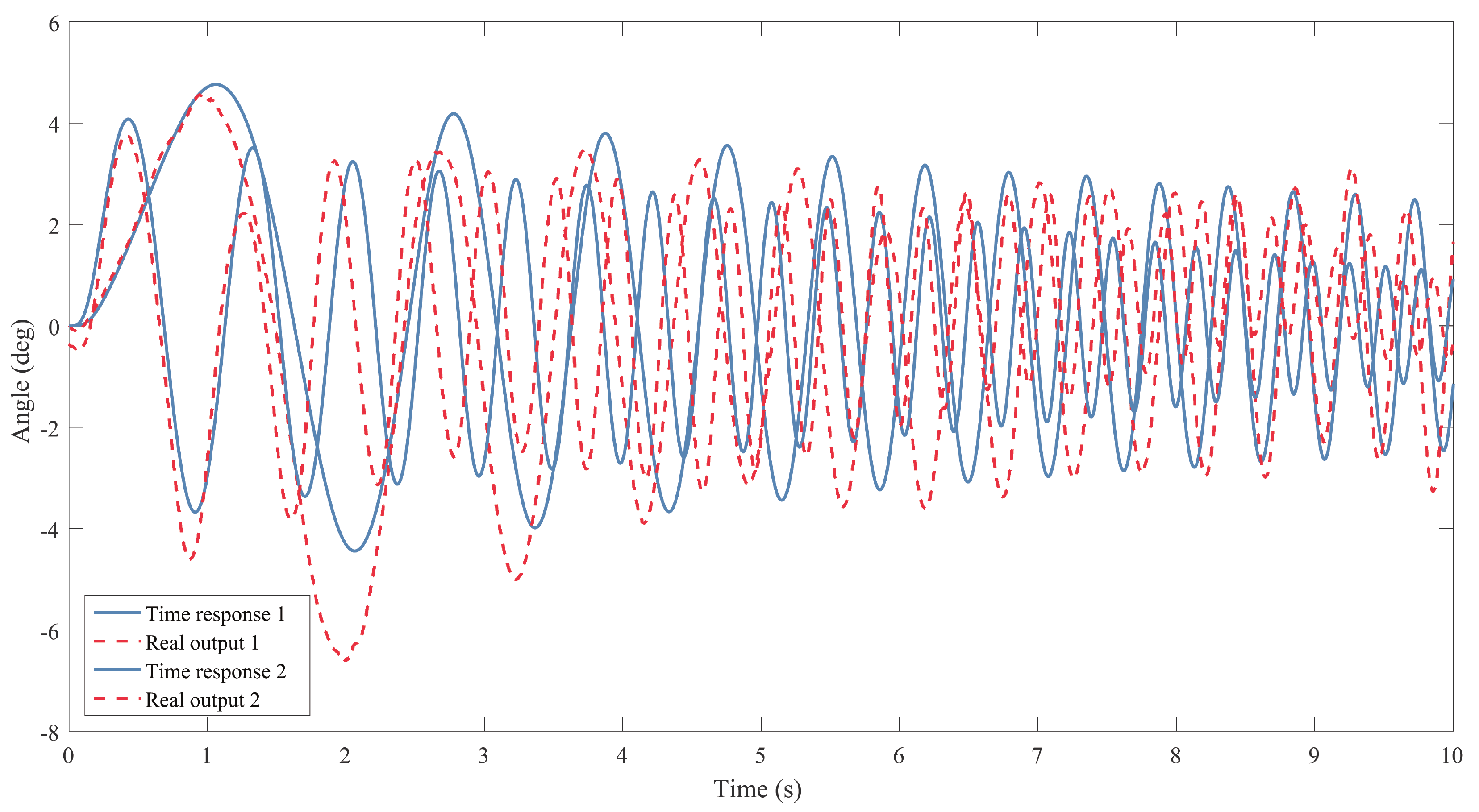}
        \caption{Roll}
        \label{fig282}
    \end{subfigure}
    \caption{Comparison of transfer functions time-response and real data (for the same input)}\label{fig28}
\end{figure}

In order to verify the nominal models, the real experimental data (dashed) and the time-domain response (solid) are depicted in Fig. \ref{fig281} and \ref{fig282}. For the sake of brevity, only two of the experiment results are shown. In these figures, the applied inputs to the nominal model are the same as the real experiments.  

\subsection{Uncertainty weighting function}
As mentioned before, based on Eq. (\ref{eq21}), $W$ is a fixed transfer function for normalizing $\Delta$. In other words, $W\Delta$ expresses the normalized system variation away from unity at each frequency, and based on Eq. (\ref{eq21}), may be obtained by $G\left(jw\right)/G_0\left(jw\right)-1=\Delta\left(jw\right) W\left(jw\right)$. Since $\lVert\Delta\rVert_\infty<1$, the following can be concluded, which provide a practical approach for obtaining the uncertainty weighting function:
\begin{equation}
\lvert\frac{G\left(jw\right)}{G_0\left(jw\right)}-1\rvert \leq\lvert W\left(jw\right)\rvert \qquad \forall w.
\end{equation}

\begin{figure}
    \centering
    \begin{subfigure}[b]{0.48\textwidth}
         \includegraphics[height=0.9\textwidth,width=\textwidth]{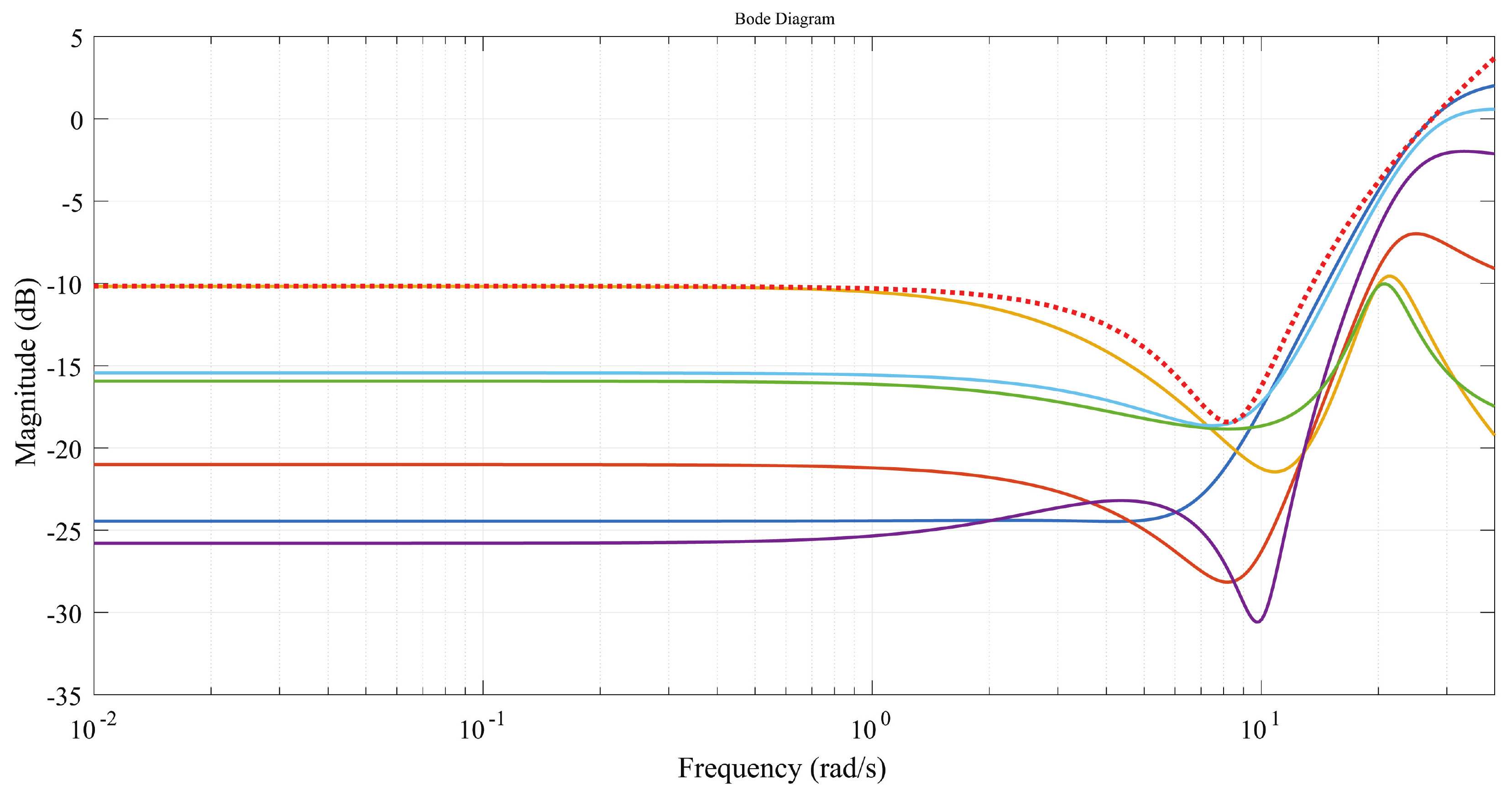}
        \caption{Pitch}
        \label{fig291}
    \end{subfigure}
    ~ 
    \begin{subfigure}[b]{0.48\textwidth}
         \includegraphics[height=0.9\textwidth,width=\textwidth]{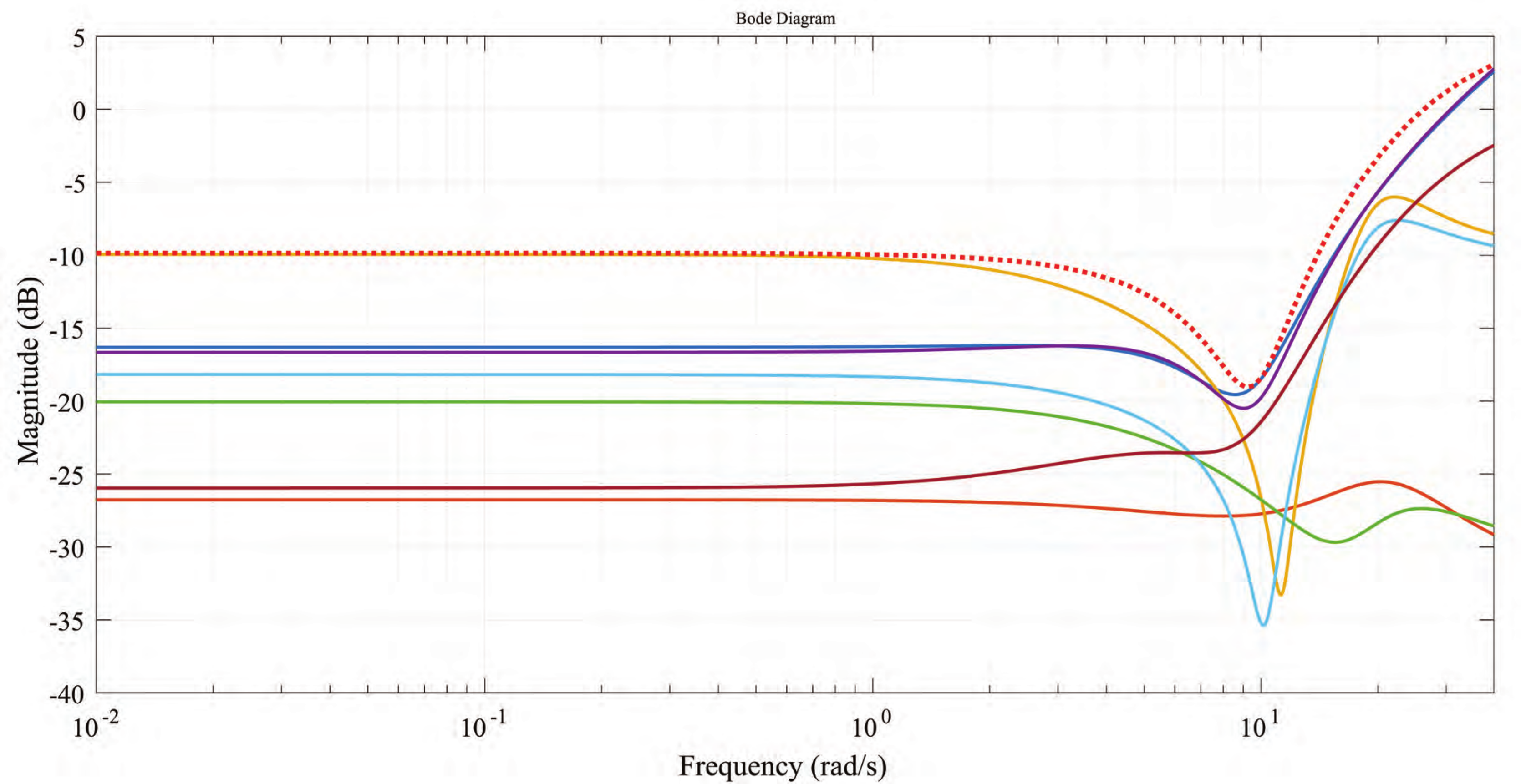}
        \caption{Roll}
        \label{fig292}
    \end{subfigure}
    \caption{Multiplicative uncertainty profiles}\label{fig29}
\end{figure}

After selecting a nominal model between all best-fitted transfer functions, other systems can be considered as a perturbed model of the nominal system, and $\lvert G\left(jw\right)/G_0\left(jw\right)-1 \rvert$ is plotted for each experiment (as illustrated in Fig. \ref{fig29}). As a result, by estimating a transfer function, which is an upper bound for those variations, the uncertainty weighting function (the dashed red line in Fig. \ref{fig291} and \ref{fig292}) is acquired. The uncertainty weighting functions for pitch and roll angles are obtained as follows:
\begin{align}
W_{pitch}=\frac{1659.6\left(s^2+2.868s+60.44\right)}{\left(s+2.477e04\right)\left(s+9.678\right)},\\
W_{roll}=\frac{1.9017\left(s^2+3.813s+91.61\right)}{s^2+43.53s+545.3}.
\end{align}

The above-mentioned steps for obtaining the nominal models of the systems and the uncertainty weighting functions are summarized as below:
\begin{enumerate}[
	leftmargin=*,
	label={\textbf{\textit{Step \arabic*.}}}
	]
	
	\item Designing appropriate experiments and obtaining 8 and 9 data packages for pitch and roll angles, respectively.
	\item Fitting 15 transfer functions for each data package.
	\item Selecting the best-fitted transfer function for each data package.
	\item Selecting the nominal models among 8 and 9 transfer functions obtained for pitch and roll angles.
	\item	Considering the other transfer functions as perturbed models, obtaining the  $\lvert G\left(jw\right)/G_0\left(jw\right)-1 \rvert$s, and plotting them for roll and pitch angles.
	\item	Obtaining the uncertainty weighting functions, which are an upper bound for the variations obtained in {\textbf{\textit{Step 5}}}, for
	each roll and pitch angle.
\end{enumerate}

\section{Robust ${H}_\infty$ Control ُSynthesis}

\begin{figure} [h]
    \centering
    \includegraphics[height=2.0in]{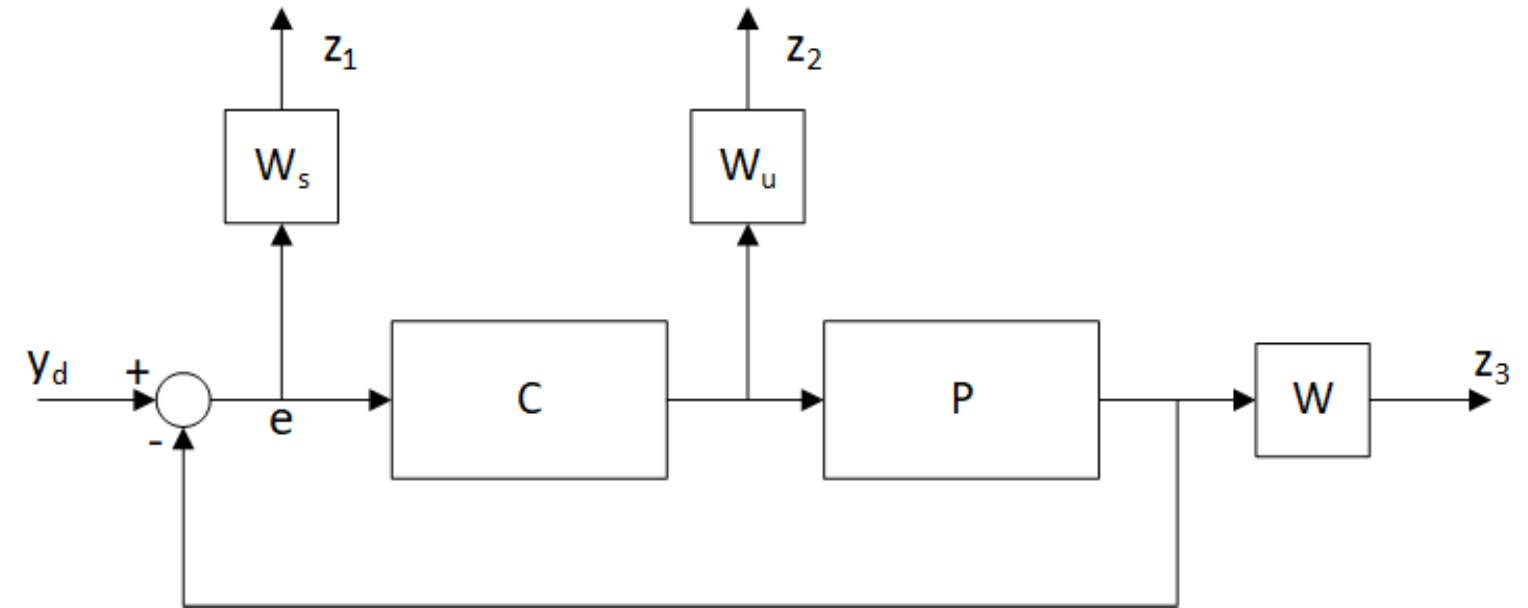}
    \caption{Block diagram of closed-loop system as a generalized regulator problem in $H_\infty$ framework}
    \label{fig31}
\end{figure}

The goal of ${H}_\infty$ controller is providing a robust stable closed-loop system with high performance tracking and disturbance rejection, in the presence of uncertainty and the actuator saturation. Based on Fig. \ref{fig31}, in order to achieve this goal the following objectives shall be met simultaneously:
\begin{enumerate}
\item $\rVert T_{z_1y_d}\lVert_\infty=\rVert W_sS\lVert_\infty\le 1$, where  $T_{z_1y_d}=\frac{W_s}{1+CP}=W_sS$ is the transfer function from $y_d$ to $z_1$, which is the nominal tracking performance in an asymptotic sense.  $W_s$ is a frequency dependent sensitivity weighting function for normalizing and shaping the closed loop performance specification. 
\item  $\rVert T_{z_2y_d}\lVert_\infty=\rVert W_uU\lVert_\infty\le 1$, in which  $T_{z_2y_d}=\frac{W_uC}{1+CP}=W_uCS=W_uU$ is the transfer function from $y_d$ to $z_2$, and expresses the nominal performance of the control effort. $W_u$ is a frequency dependent weighting function for normalizing and  shaping the control input.
\item $\rVert T_{z_3y_d}\lVert_\infty=\rVert WT\lVert_\infty\le 1$, which is the result of the small gain theorem. $T_{z_3y_d}=\frac{WCP}{1+CP}=WT$ is the transfer function from $y_d$ to $z_3$, and $W$ represents the uncertainty weighting function. 
\end{enumerate}
Considering this system with  input, $y_d$, and  output vector $\mathbf{z}=[z_1, z_2, z_3]^T$, the aforementioned conditions can be merged in an induced norm of the transfer function from  $y_d$ to  $\mathbf{z}$, in which the goal is to find a controller to minimize this norm. This problem is called a mixed-sensitivity problem and is formulated as follows:
\begin{equation}
\gamma_{opt}= \min_{C}{\left\lVert T_zy_d \right\rVert}=\min_{C}{\left\lVert\begin{matrix}
           W_sS \\
           W_uU \\
          WT
         \end{matrix}\right\rVert_\infty} .
\end{equation}
In order to define the sensitivity weighting function, first an ideal closed-loop transfer function is designed. In this case, based on the standard second order systems, $\omega_n$ and $\zeta$ are chosen such that the following overshoot and the settling time are satisfied
\begin{align}
&T_{cl}=\frac{\omega_n^2}{s^2+2\zeta \omega_ns+\omega_n^2},\\
&t_s\approx\frac{4}{\zeta\omega_n}=0.3,\\
&M_p=e^{\frac{-\pi\zeta} {\sqrt{1-\zeta^2}}}=10^{-6}, \qquad 0<\zeta<1 .
\end{align}
To achieve $W_s$, we use the inverse of the desired sensitivity function:
\begin{align}
&T_{id}=\frac{247.3}{s^2+30.67s+247.3}\label{tid},\\
&W_s=a\frac{1}{S_{id}}=a\frac{1}{1-T_{id}}=a\frac{s^2+30.67s+247.3}{s\left(s+30.67\right)} \label{eq317},
\end{align}
where $a\le1$ is a tuning parameter. Eq. (\ref{eq317}) is modified by inserting a non-dominant pole and by slightly transferring its purely imaginary pole to the left half plane in order to convert $W_s$ to a strictly proper and stable weighting function, respectively. Hence, $W_s$ is corrected as:
\begin{equation}\label{eq318}
W_s=a\frac{\left(s^2+30.67s+247.3\right)}{\left(s+1000\right)\left(s+30.67\right)\left(s+0.001\right)}.
\end{equation}

In the obtained nominal system, Eq. (\ref{eq27}), the input of the system is the desired angle, which should be less than $20^\circ$. Therefore, $W_u$ is set to $0.05$. Having set the weighting functions, the mixed sensitivity problem can be solved using the robust toolbox of MATLAB. After a few correction steps, the controller for the pitch angle by tuning $a=0.88$ with  $\gamma_{opt}=0.9929$ is synthesized as follows:
\begin{equation*}\resizebox{0.98\hsize}{!}{$
C_{pitch}\left(s\right)=\frac{3.3227e05\left(s+ 2.477e04\right)\left(s+ 461.2\right)\left(s +50.77\right)\left(s+ 9.678\right)\left(s +5.373\right)\left(s^2+10.12s+ 390.4\right)}
{\left(s+2.477e04\right) \left(s+1.673e04\right) \left(s+1000 \right)\left(s+30.67\right) \left(s+8.082\right) \left(s+0.001\right) \left(s^2+80.13s + 3556\right)}.
$}\end{equation*}

\begin{figure} [h]
    \centering
    \includegraphics[width=\textwidth]{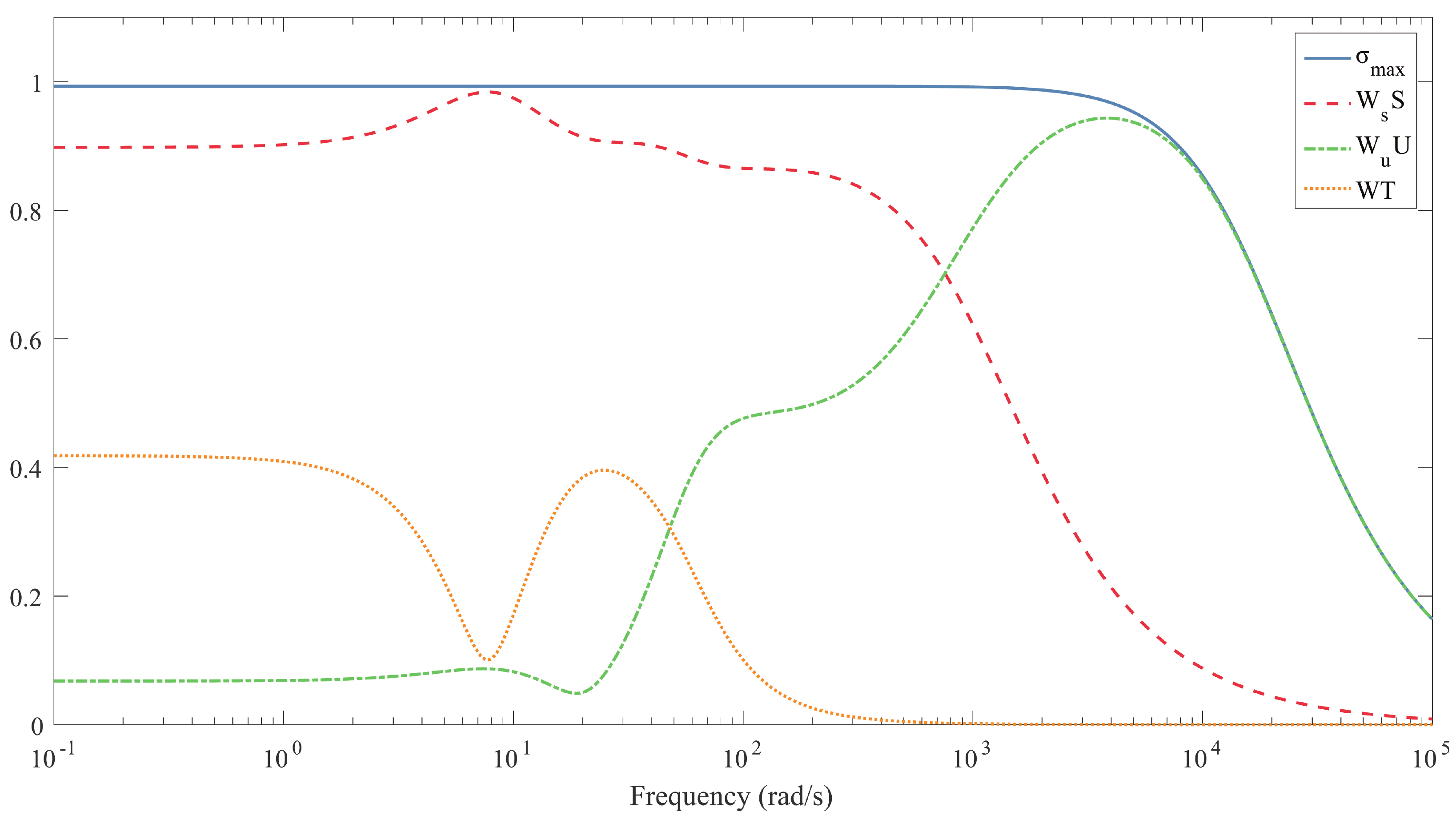}
    \caption{The closed loop system singular values and Bode plot of $W_sS$, $ W_uU$, $WT$ (Pitch)}
    \label{fig32}
\end{figure}

The singular values of the closed-loop system,  and the frequency response of each I/O pair are depicted in Fig. \ref{fig32}. The solid blue line shows the maximum singular value of the closed loop system, which is flat within a large frequency range, and is less than one in all range of the frequency domain.   The dotted orange line shows the Bode diagram of $WT$ in which the maximum value is 0.42  indicating the robustness of the closed loop system with a margin of greater than two. The performance transfer function and the control effort transfer function are shown by the red and green lines, respectively. Based on  Fig. \ref{fig32}, the magnitude of $WT$  starts to reduce significantly in frequencies higher than about 30 Hz due to the limited bandwidth of the system. About this frequency, a slight decrease also appears in the magnitude of $W_sS$ until roughly 900 Hz where a dramatic drop occurs. Conversely, the magnitude of $ W_uU$ increases in these frequencies to compensate for the lack of the system stability and performance. However, since the robot working frequency is less than approximately 20 Hz, high frequencies performance of the system can be ignored. In low frequencies, as expected, the magnitude of $W_sS$  increases when that of $WT$ decreases.

\begin{figure} [h]
    \centering
    \includegraphics[width=\textwidth]{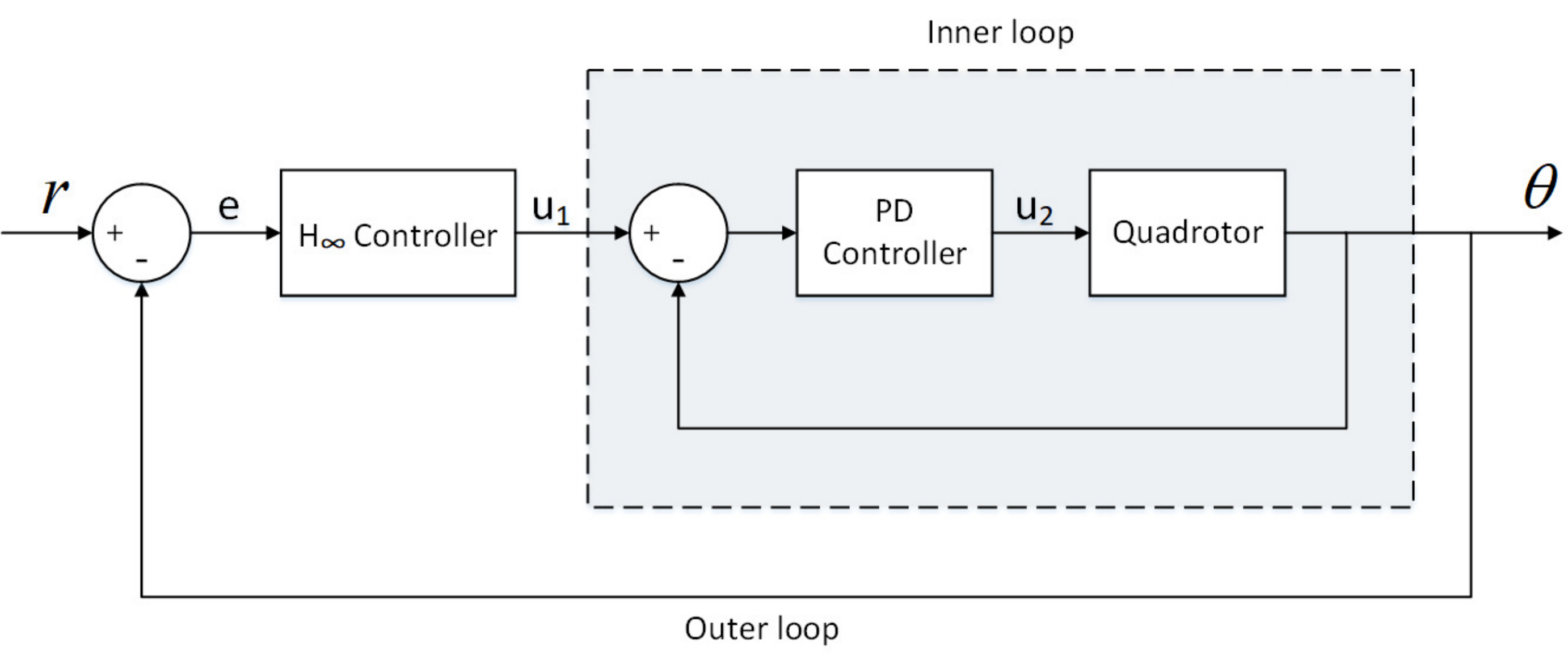}
    \caption{Block diagram of the control structure}
    \label{fig36}
\end{figure}

\begin{figure}
    \centering
    \begin{subfigure}[b]{0.48\textwidth}
         \includegraphics[height=0.8\textwidth,width=\textwidth]{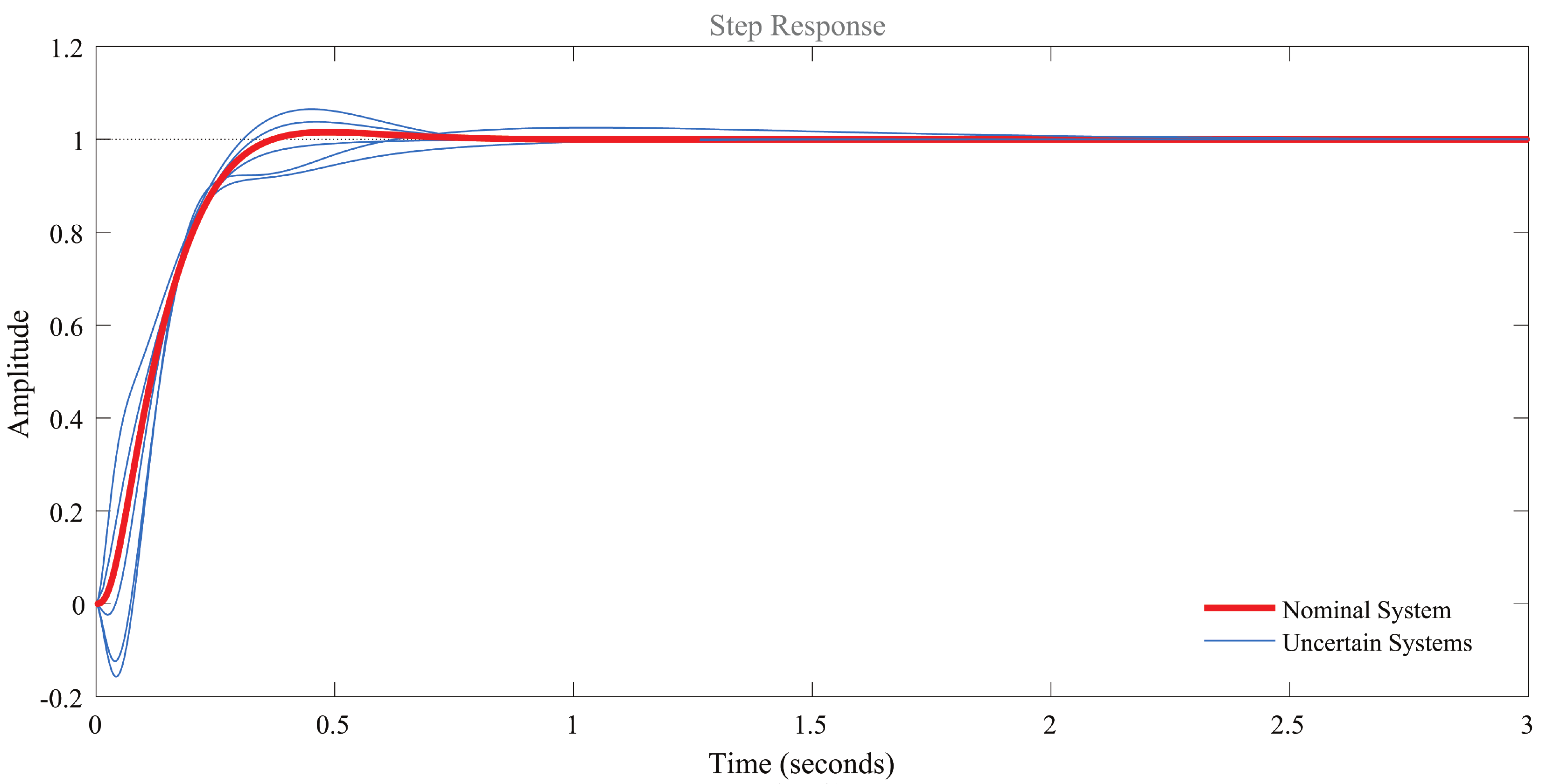}
        \caption{Time response}
        \label{fig331}
    \end{subfigure}
    ~ 
    \begin{subfigure}[b]{0.48\textwidth}
         \includegraphics[height=0.76\textwidth,width=\textwidth]{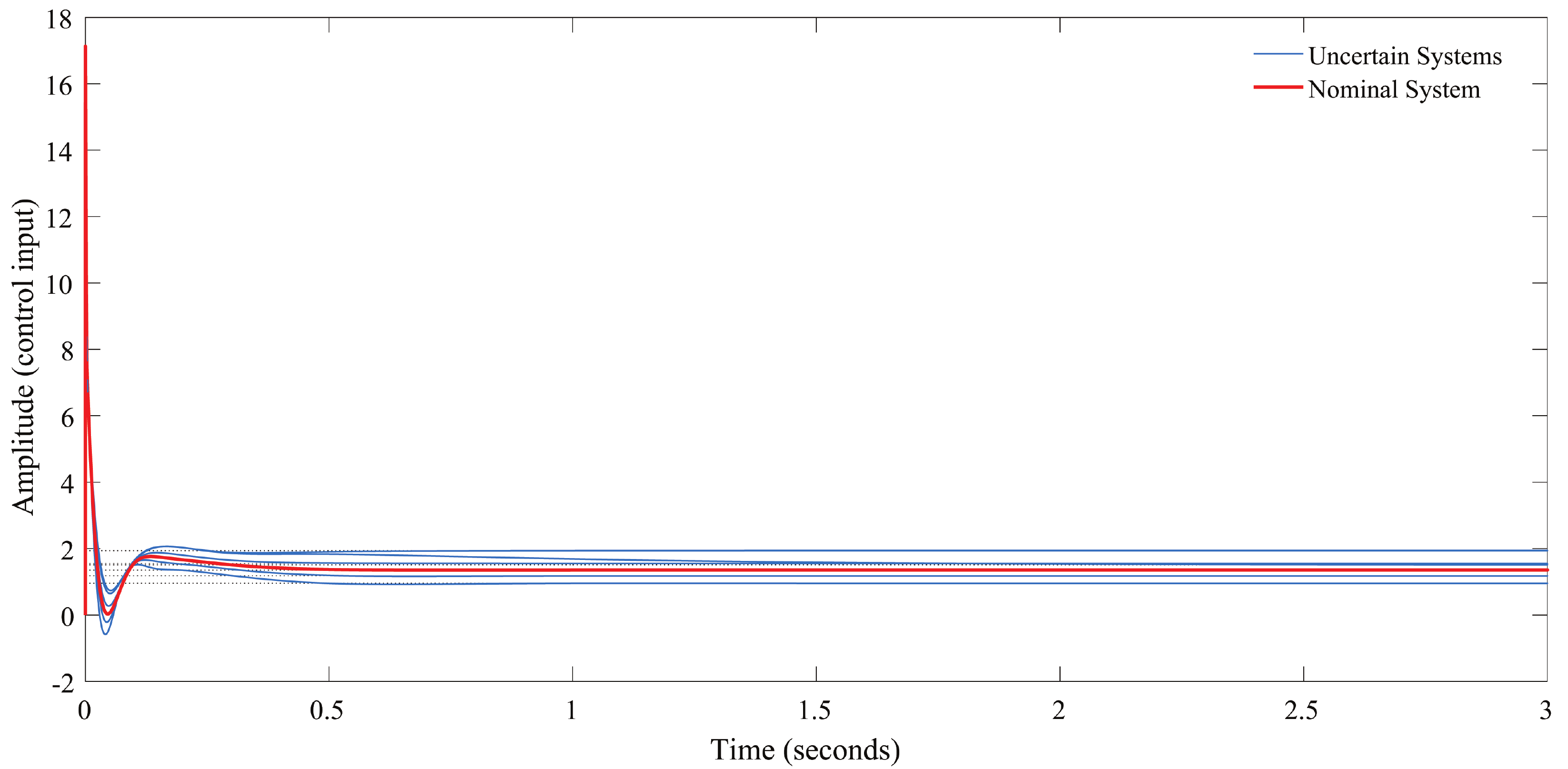}
        \caption{Controll effort}
        \label{fig332}
    \end{subfigure}
    \caption{Closed-loop control of nominal and uncertain systems for the unit-step input (Pitch)}\label{fig33}
\end{figure}

As mentioned in section \ref{sec2}, a PD controller has been employed to identify the system. Therefore, to implement the final controller a cascade architecture must be used in which the PD controller is in the inner loop, and the $H_\infty$ controller is functioning in the outer loop as illustrated in  Fig. \ref{fig36}.
The closed-loop step response of the nominal (red) and some uncertain samples of the system are shown in Fig. \ref{fig331}. It can be seen that the controller can robustly stabilize the system and provide a fast response with a reasonable control effort, which is depicted in Fig. \ref{fig332}. 

\begin{figure} [h]
    \centering
    \includegraphics[width=\textwidth]{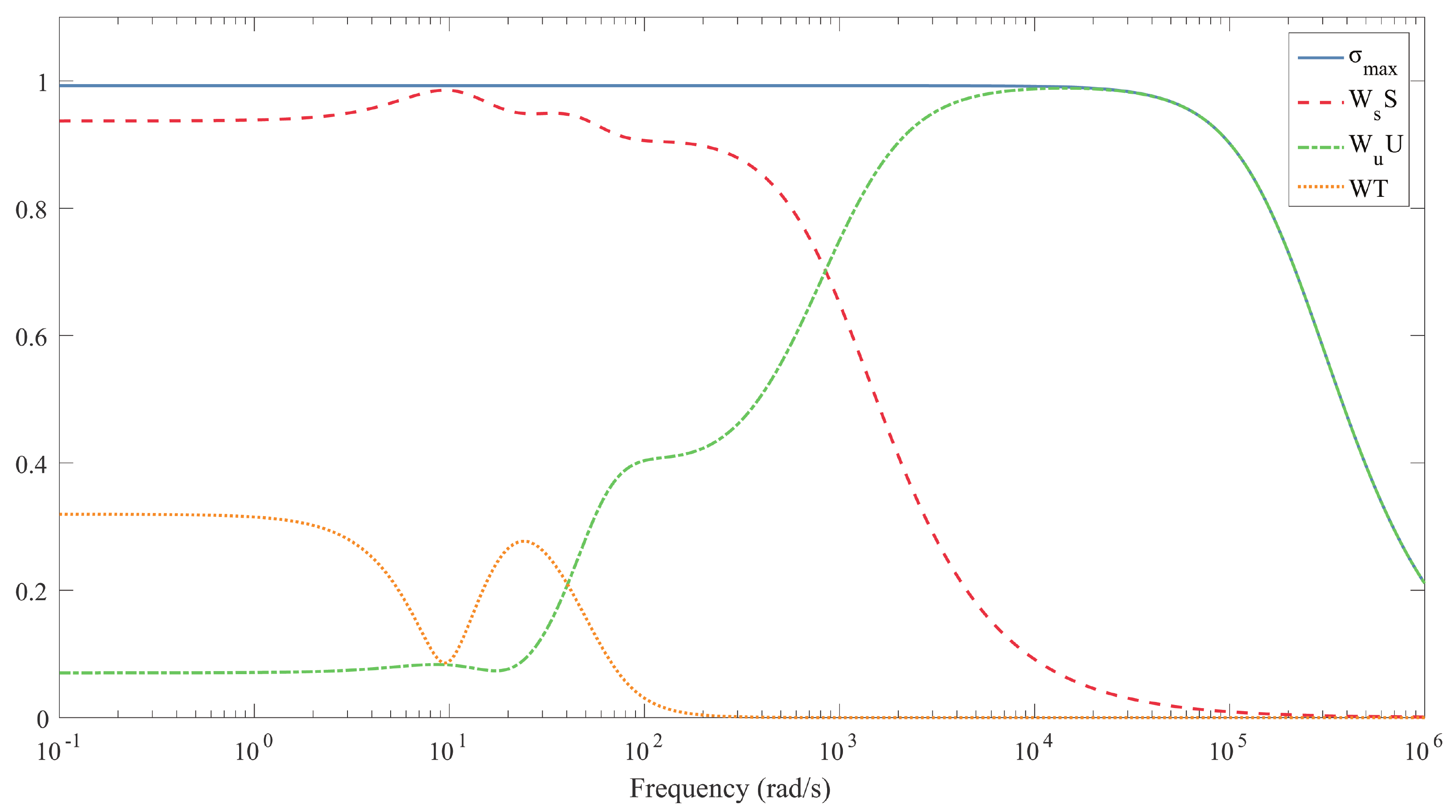}
    \caption{The closed loop system singular values and Bode plot of $W_sS$, $ W_uU$, $WT$ (Roll)}
    \label{fig34}
\end{figure}

\begin{figure}
    \centering
    \begin{subfigure}[b]{0.48\textwidth}
         \includegraphics[height=0.8\textwidth,width=\textwidth]{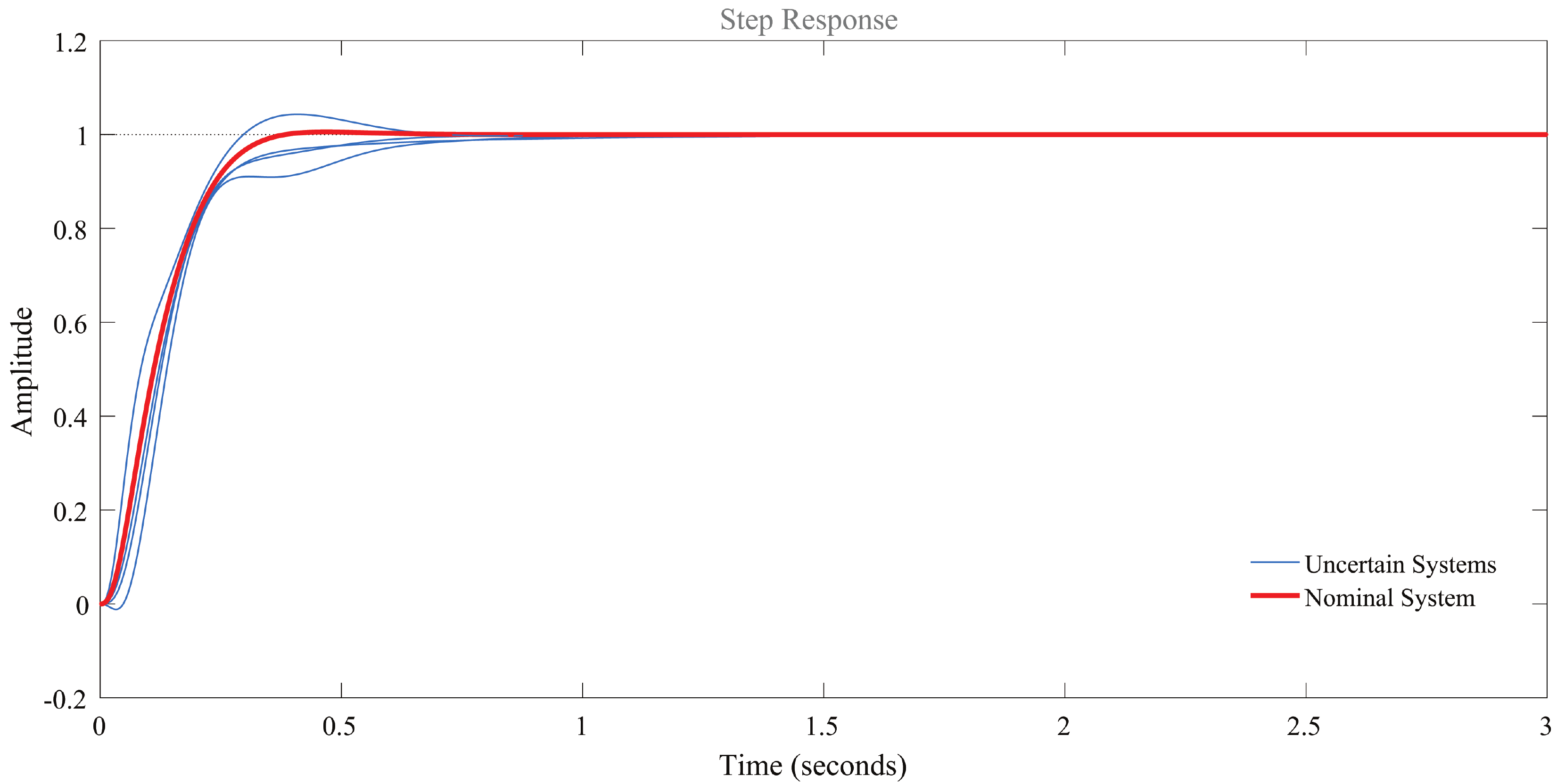}
        \caption{Time response}
        \label{fig351}
    \end{subfigure}
    ~ 
    \begin{subfigure}[b]{0.48\textwidth}
         \includegraphics[height=0.77\textwidth,width=\textwidth]{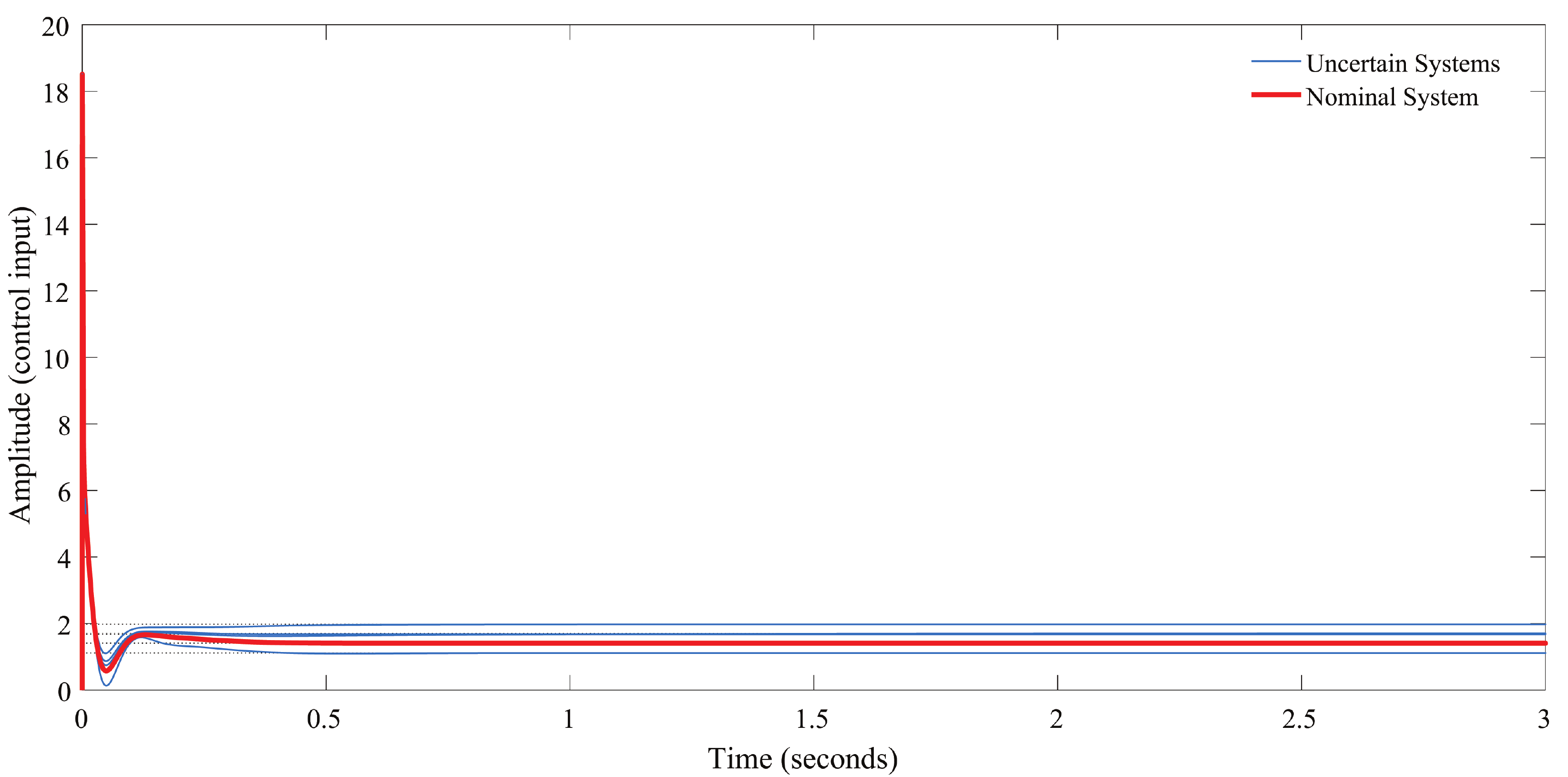}
        \caption{Controll effort}
        \label{fig352}
    \end{subfigure}
    \caption{Closed-loop control of nominal and uncertain systems for the unit-step input (Roll)}\label{fig35}
\end{figure}

The same procedure is applied to obtain the $H_\infty$ controller for the roll angle.  $W_u$ is set to 0.05, and $W_s$ is same as that of the pitch angle, in Eq. (\ref{eq318}). After solving the mixed-sensitivity problem,  $\gamma_{opt}$ value is obtained 0.9925, for $a = 0.92$ and the following controller:
\begin{equation*}\resizebox{0.98\hsize}{!}{$
C_{roll}\left(s\right)=\frac {4.3173e06 \left(s+371.5\right) \left(s+59.89\right) \left(s+6.764\right) \left(s^2 + 43.53s + 545.3\right)\left(s^2 + 19.03s + 426.2\right)}
{\left(s+2.175e05\right) \left(s+1000\right) \left(s+30.67\right) \left(s+27.95\right) \left(s+19.41\right) \left(s+0.001\right)\left(s^2 + 80.99s + 3749\right)} .
$}\end{equation*}
The singular values of the closed-loop system and the frequency response of each I/O pair are illustrated in Fig. \ref{fig34}. The analysis of this plot is the same as that of the pitch angle, which is omitted for the sake of brevity. The step responses of the nominal and perturbed systems in the presence of $H_\infty$ controller are depicted in Fig. \ref{fig351}. The control efforts  are also shown in Fig. \ref{fig352}. 

\section{Experimental results}
\subsection{Real Implementation}
In this section, we verify the performance and the applicability of the controllers which were designed in the preceding section by implementing on the real robot. 
The control structure is the same as Fig. \ref{fig36} for both pitch and roll angles. By setting the time-step of the system to 0.001, the continuous controllers are discretized with the frequency of $1\, kHz$ for the practical implementation. To make this task easier, an order reduction is performed that decreases the order of the controllers from 8 to 6.  Finally, the controllers for both axes after applying the order reduction are given by:
\begin{align*}
&C_{pitch}=\frac{3.3227e05 (s+461.5) (s+45.87) (s+5.903) (s^2 + 9.836s + 388.1)}{(s+1.673e04) (s+1000) (s+25.13) (s+0.001)(s^2 + 79.66s + 3577)},\\
&C_{roll}=\frac{4.3173e06 (s+371.8) (s+53.33) (s+7.035) (s^2 + 18.34s + 428.1)}{(s+2.175e05) (s+1000) (s+27.87) (s+0.001) (s^2 + 80.84s + 3820)}.
\end{align*}
Based on the frequency response of each I/O pair, this reduction has a negligible effect.

Having implemented the controllers, some tests are performed in the outdoor and indoor environments to evaluate the $H_{\infty}$ controller ability in coping with uncertainties. In the first experiment, the  operator drives the quadrotor using a radio-controller to assess the qualitative performance of the robot in comparison with a well-tuned PID controller. The tuning parameters of the PID controller for both roll and pitch angles are selected as follows:
\begin{equation*}
K_p=2.6, K_i=0.2, K_d=0.65
\end{equation*}

 According to this experience, the maneuverability has been enhanced, and the quadrotor can quickly regulate the roll and the pitch angle in the presence of the external disturbances such as the unexpected wind effects (i.e. caused by the drone itself while flying near the walls or the weather) or some manual disturbances applied by the operator. Since the states of the quadrotor are coupled, a refinement on the altitude control of the robot is also observed. The video of this experiment is also available online\footnote{\href{https://youtu.be/vAi4x2_XSQQ}{https://youtu.be/vAi4x2\_XSQQ}}.

\begin{figure}
    \centering
    \begin{subfigure}[b]{0.48\textwidth}
         \includegraphics[height=0.9\textwidth,width=\textwidth]{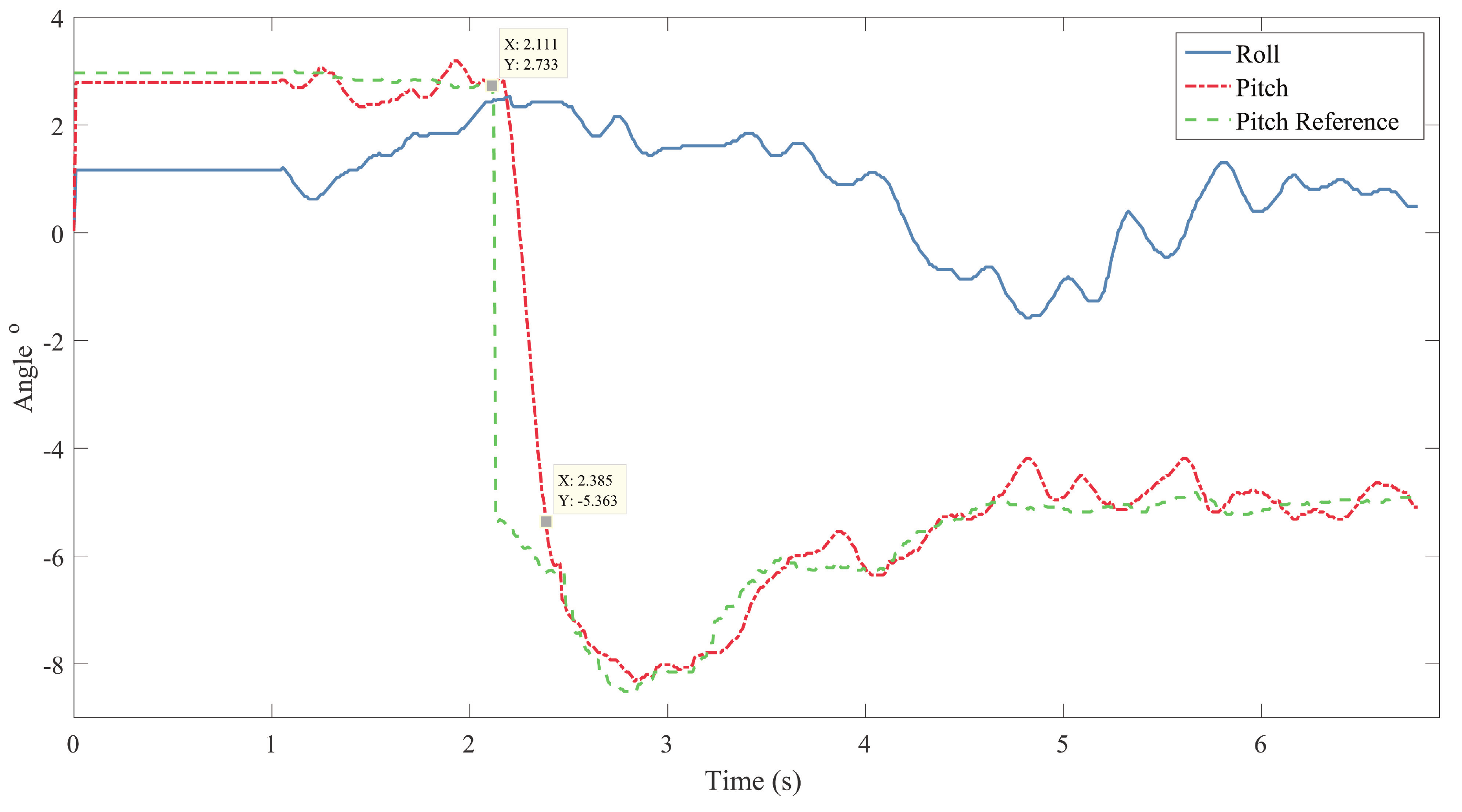}
        \caption{Robust $H_{\infty}$ controller}
        \label{fig411}
    \end{subfigure}
    ~ 
    \begin{subfigure}[b]{0.48\textwidth}
         \includegraphics[height=0.9\textwidth,width=\textwidth]{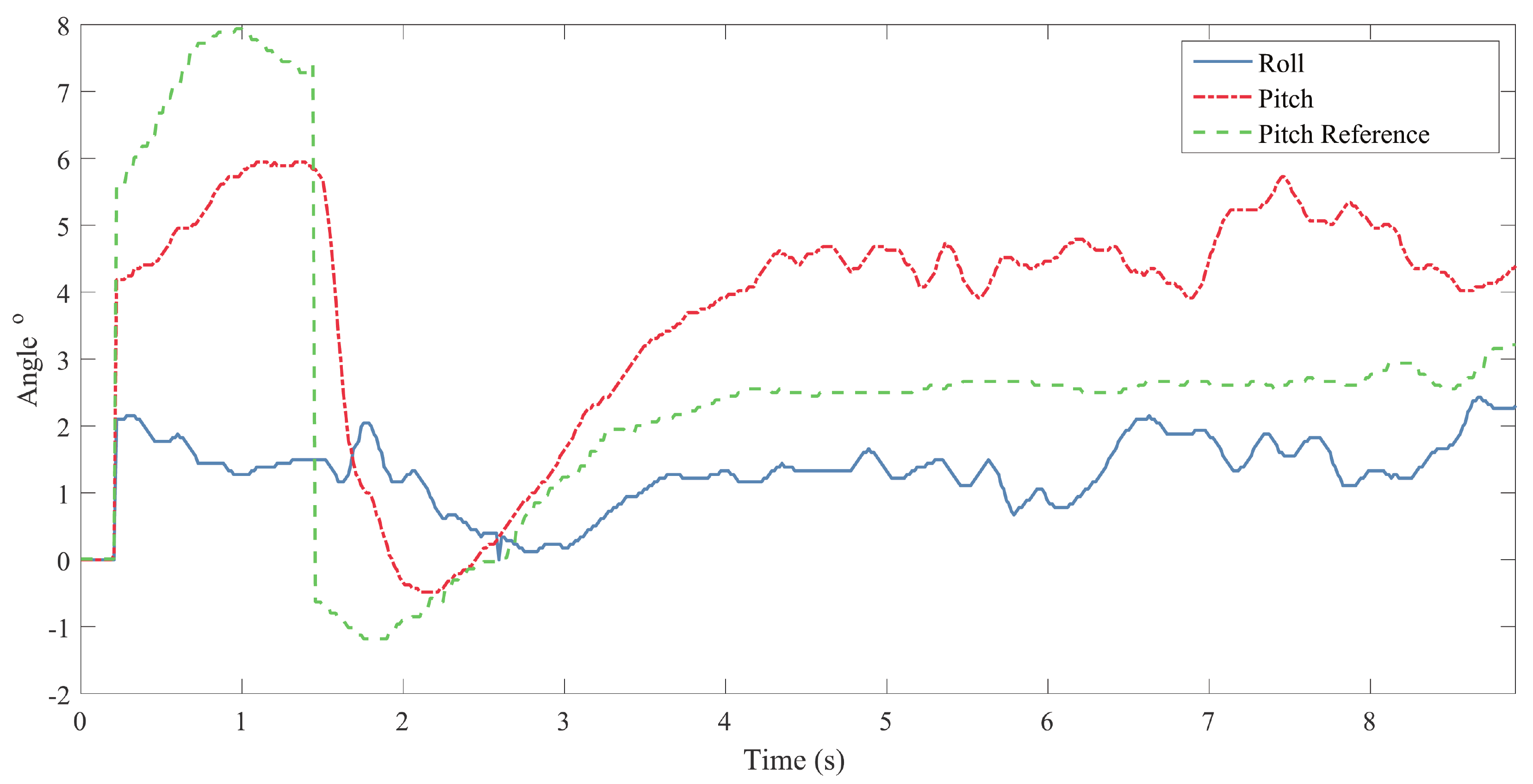}
        \caption{PID controller}
        \label{fig412}
    \end{subfigure}
    \caption{Tracking and regualting performance of $H_{\infty}$ and PID controllers }\label{fig41}
\end{figure}

To evaluate the robustness, stability and tracking performance of the proposed controllers, an experiment is done in which the quadrotor should track the reference input in the outdoor environment, in the presence of light wind. Fig. \ref{fig411} shows the tracking performance of the pitch controller and the ability of the roll controller in regulating the roll angle close to zero. In addition, as it is shown by data-tip in Fig. \ref{fig411}, the performance of the system (e.g. settling time and overshoot)  is close to that of the desired transfer function obtained in Eq. (\ref{tid}). The performance of the well-tuned PID controller is also depicted in Fig \ref{fig412}. According to Fig. \ref{fig41}, the $H_{\infty}$ controller has significantly enhanced the performance of the system.

\subsection{Simulation Results}
After validating the practicability and efficiency of the proposed controllers through real experiments, in what follows, the performance of the controllers are evaluated in the simulation and is compared with a well-tuned PID and  another robust controller obtained by $\mu$-synthesis method. The robust stability and robust performance of the controllers are also analyzed using the structural singular values. The  gains of PID controller are as below:
\begin{equation*}
\begin{aligned}
&Roll: \;  &&K_p=1.18, K_i=10.6, K_d=0.0329\\ 
&Pitch: \;  &&K_p=1.01, K_i=10.2, K_d=0.0132 
\end{aligned}
\end{equation*}

\begin{figure}
    \centering
    \begin{subfigure}[b]{0.48\textwidth}
         \includegraphics[height=0.9\textwidth,width=\textwidth]{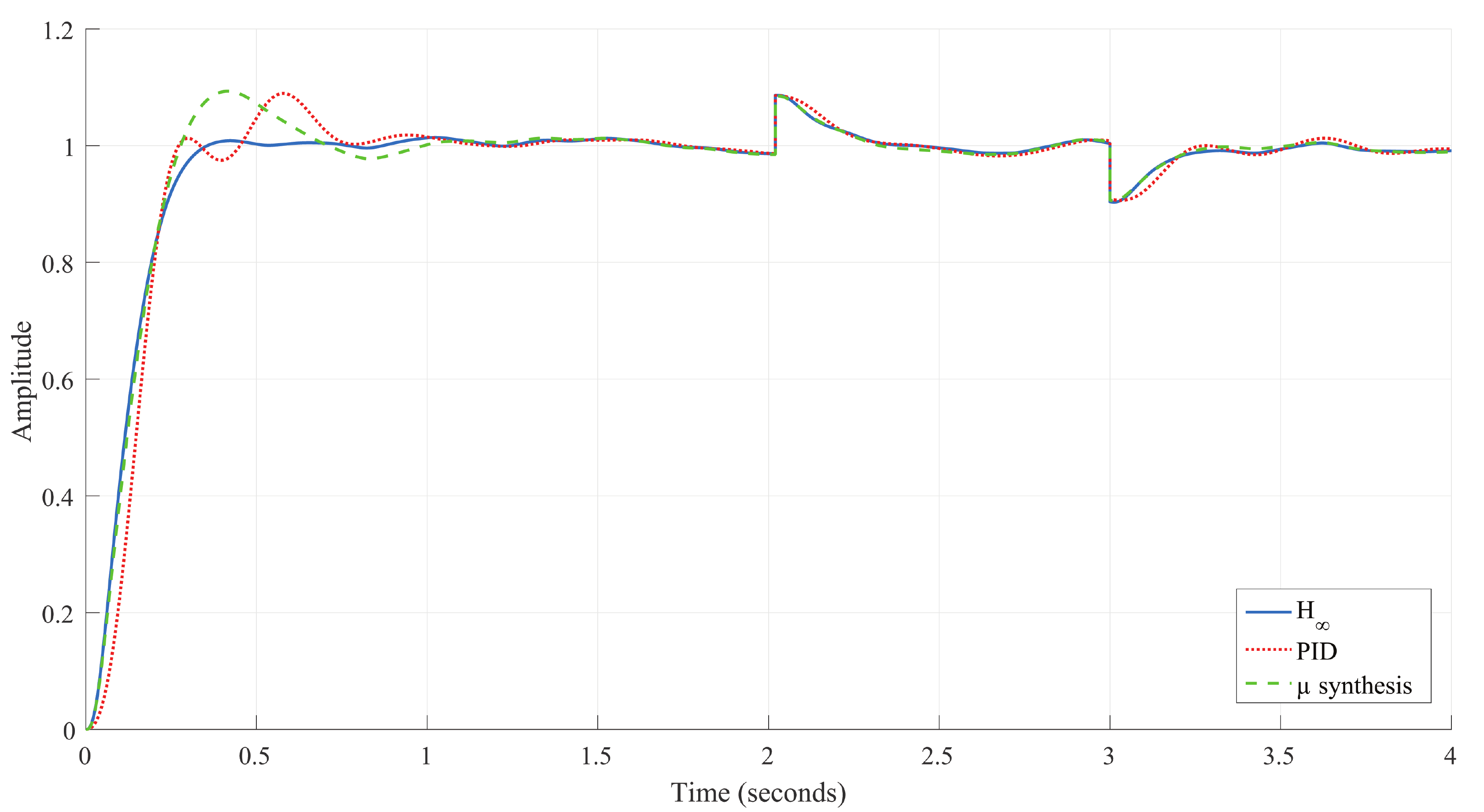}
        \caption{Pitch}
        \label{fig421}
    \end{subfigure}
    ~ 
    \begin{subfigure}[b]{0.48\textwidth}
         \includegraphics[height=0.9\textwidth,width=\textwidth]{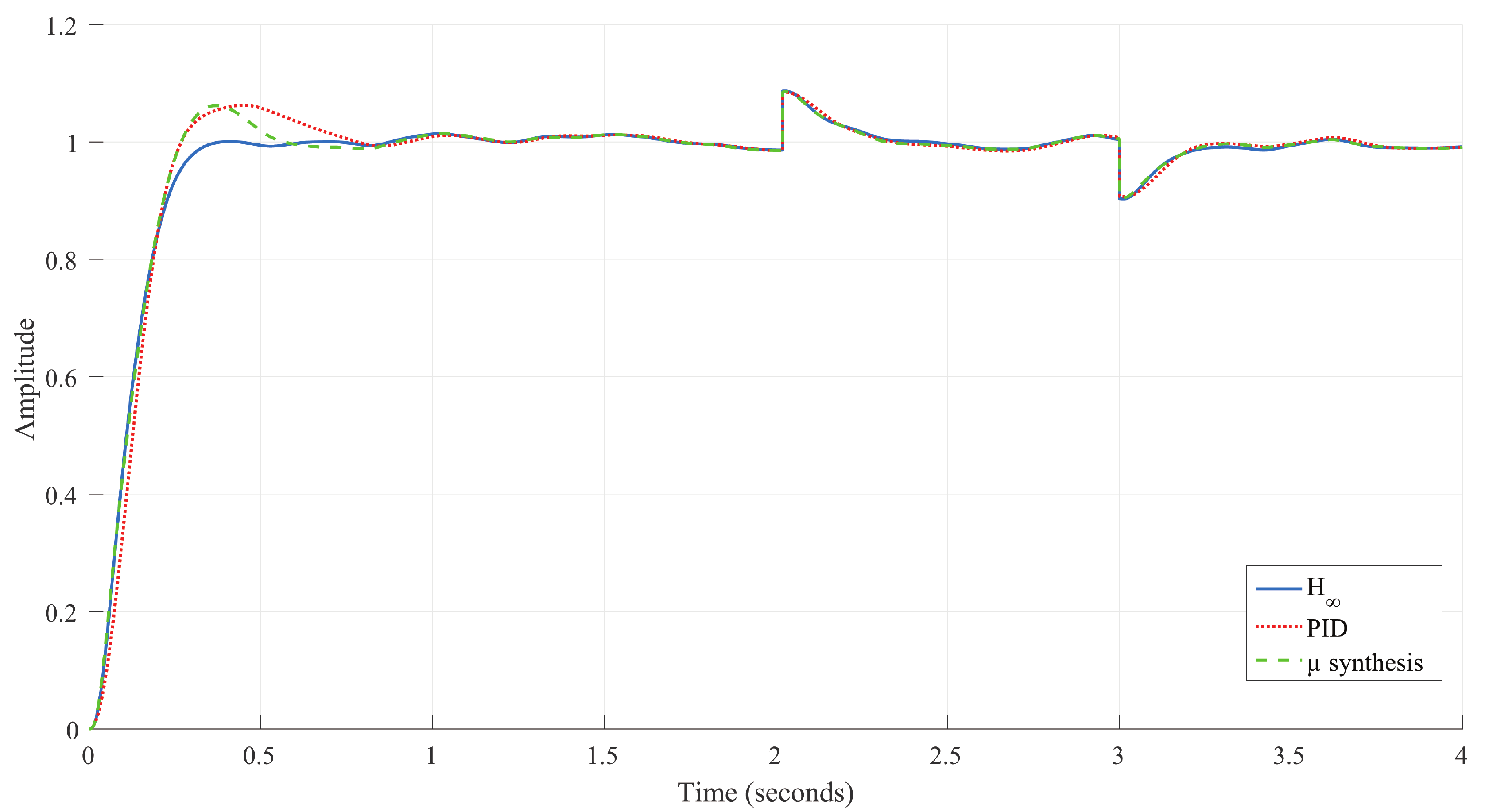}
        \caption{Roll}
        \label{fig422}
    \end{subfigure}
    \caption{$H_{\infty}$, PID and $\mu$ synthesis controllers performance in the presence of disturbance and sensor noise}\label{fig42}
\end{figure}

Fig. \ref{fig421} illustrates the step-response of an uncertain system  in the presence of the disturbances and uniform random sensor noise. The amplitude of the disturbances is $0.1$, and the maximum value of the random sensor noise is $0.02$. Fig. \ref{fig422}, depicts the roll angle performance in a similar experiment. Based on these figures, the $H_{\infty}$ controller's overshoots for the pitch and roll systems are about $1\%$ and $0.1\%$, respectively, which are less than that of PID and $\mu$ synthesis controllers which are approximately $10\%$. In addition, the $H_{\infty}$ controller regulates the system faster, and its settling time is roughly $0.35\, sec$. The settling times of both  PID and $\mu$ synthesis controllers  are about $1\,sec$.   $H_{\infty}$ controller also has better performance when disturbances occur.

\begin{figure}
    \centering
    \begin{subfigure}[b]{0.48\textwidth}
         \includegraphics[height=0.9\textwidth,width=\textwidth]{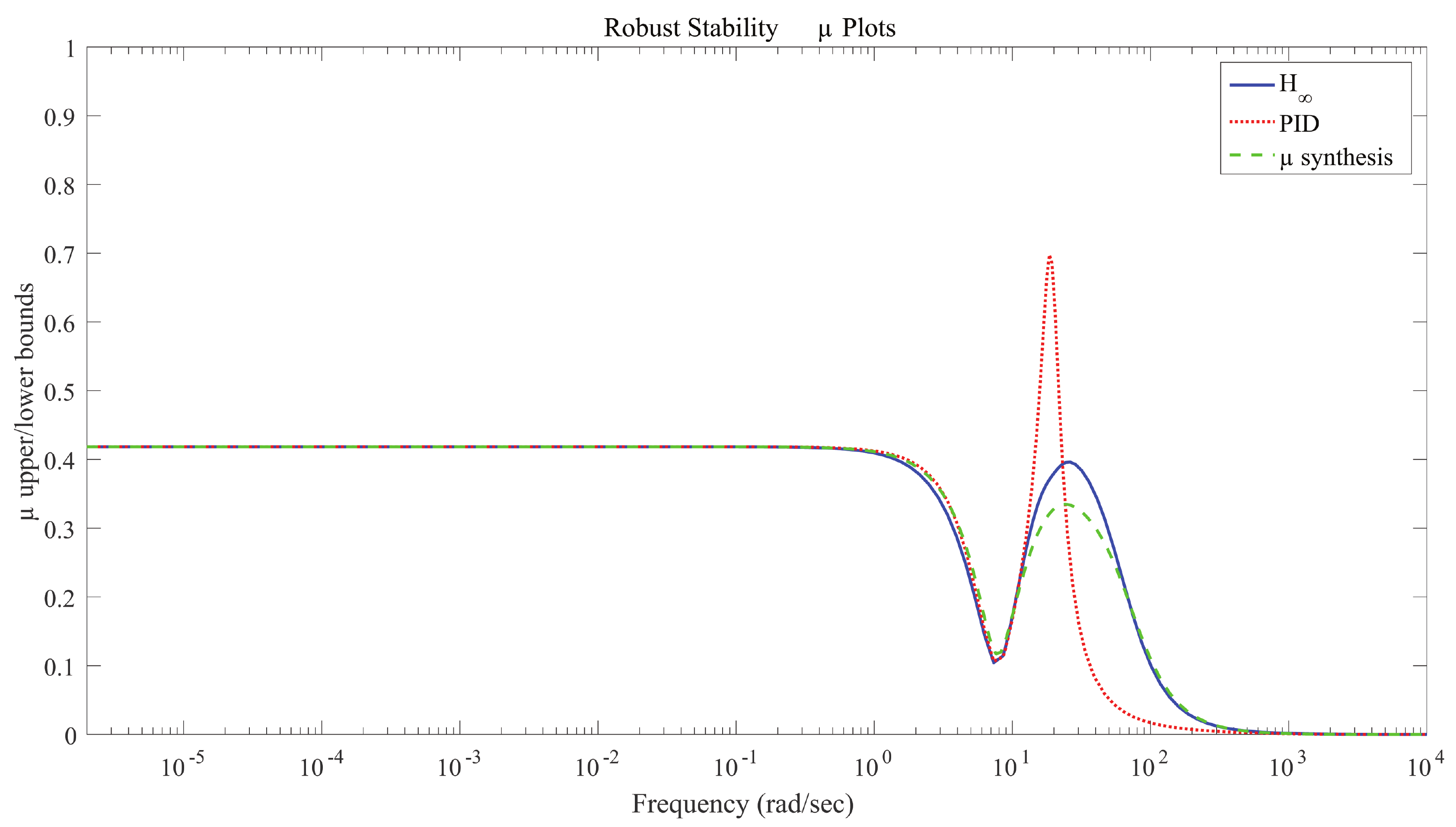}
        \caption{Pitch}
        \label{fig431}
    \end{subfigure}
    ~ 
    \begin{subfigure}[b]{0.48\textwidth}
         \includegraphics[height=0.9\textwidth,width=\textwidth]{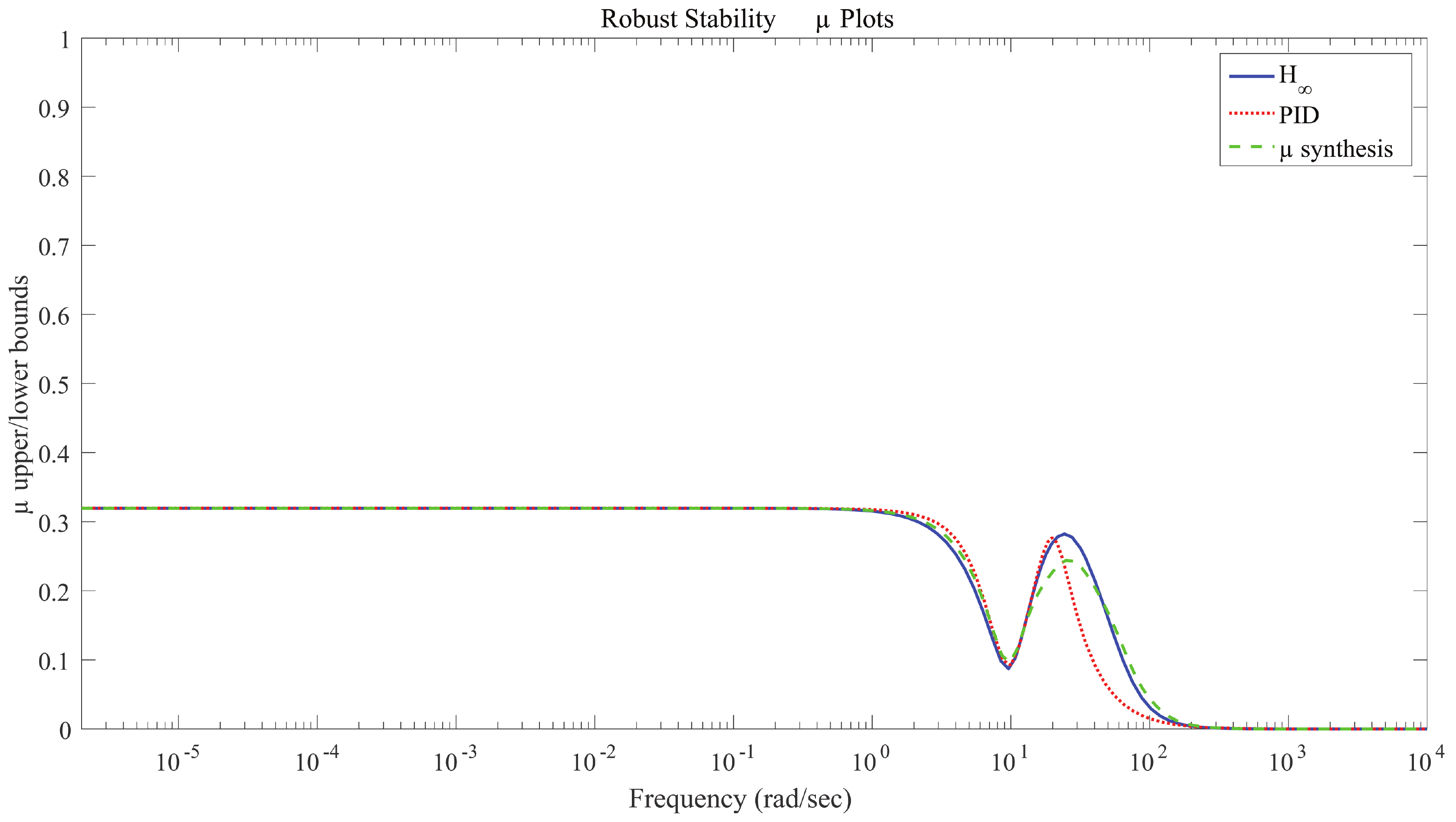}
        \caption{Roll}
        \label{fig432}
    \end{subfigure}
    \caption{Robust stability of $H_{\infty}$, PID and $\mu$ synthesis controllers }\label{fig43}
\end{figure}

\begin{figure}
    \centering
    \begin{subfigure}[b]{0.48\textwidth}
         \includegraphics[height=0.9\textwidth,width=\textwidth]{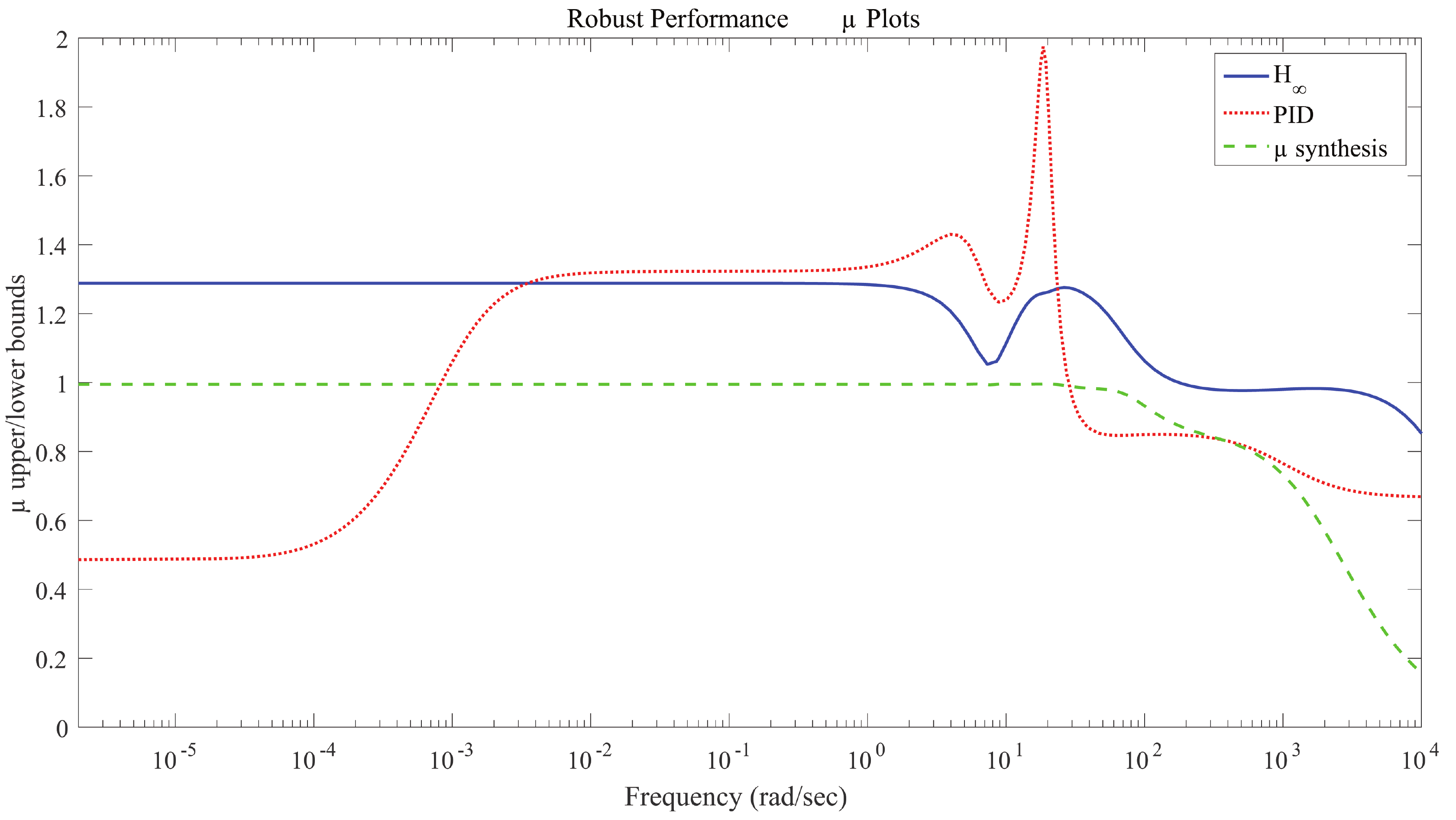}
        \caption{Pitch}
        \label{fig441}
    \end{subfigure}
    ~ 
    \begin{subfigure}[b]{0.48\textwidth}
         \includegraphics[height=0.9\textwidth,width=\textwidth]{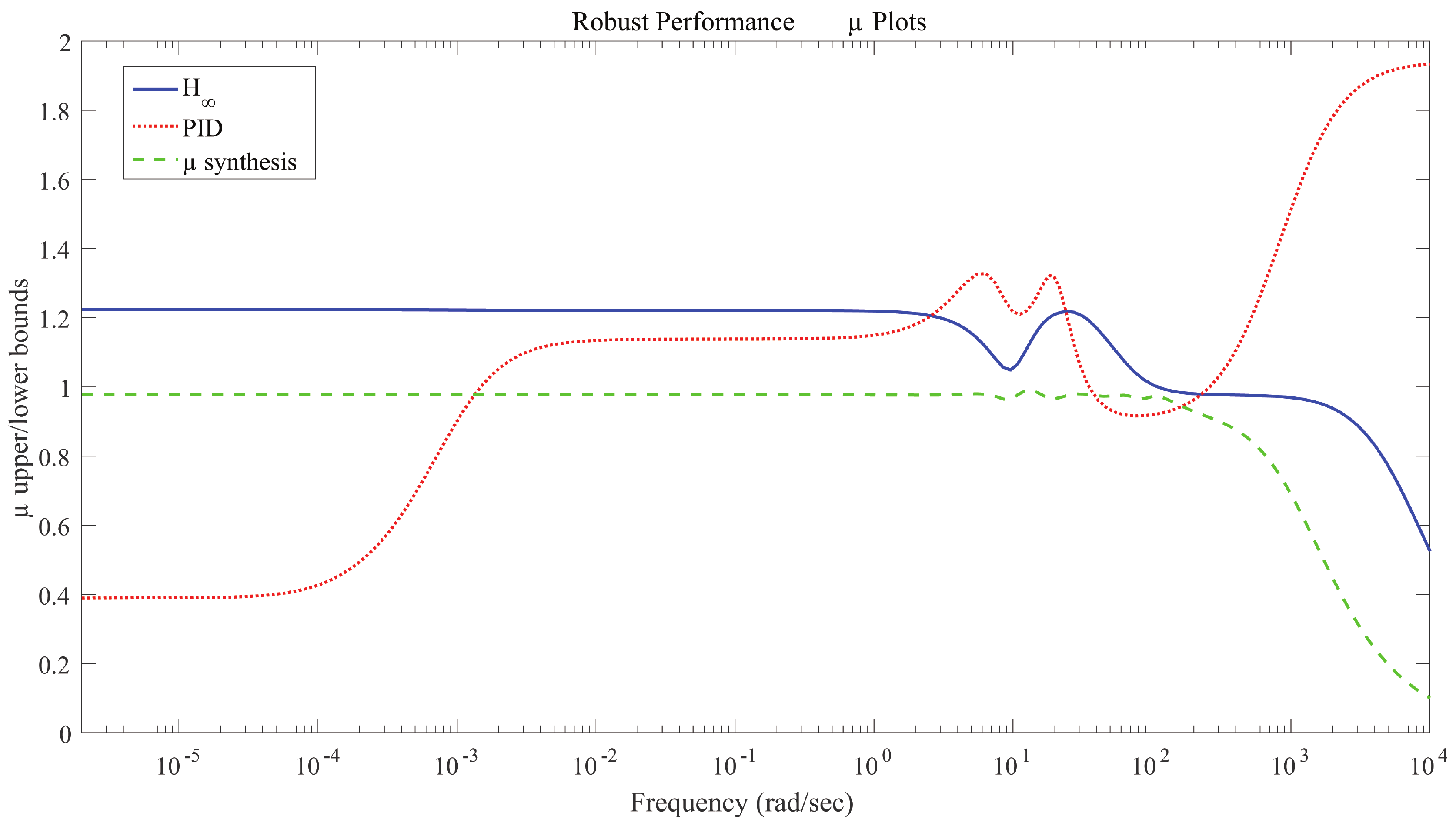}
        \caption{Roll}
        \label{fig442}
    \end{subfigure}
    \caption{Robust performance of $H_{\infty}$, PID and $\mu$ synthesis controllers }\label{fig44}
\end{figure}

In Fig. \ref{fig43}, the robust stability of the  $H_{\infty}$, $\mu$ synthesis, and PID controllers are evaluated by structural singular values. As expected, the  $H_{\infty}$ and $\mu$ synthesis controllers fulfill the robust stability condition and their $\mu$ values are less than $1$. The robust stability condition is also met for the PID controller, however, in frequencies about $20 \, rad/sec$ gets close to $1$. Moreover, when the PID parameters are set in such a way that the system has a faster response, the robust stability will be lost.
The robust performance analysis of the roll and pitch systems using structural singular values are depicted in Fig. \ref{fig44}. Predictably, the $\mu$ synthesis controller offers the robust performance in contrast to $H_{\infty}$ and PID controllers. The $\mu$ synthesis controller is a conservative approach to have the robust performance, which results in poor time-domain performance in comparison with the $H_{\infty}$ controller. Moreover, $\mu$ synthesis method produces a high-order controller which makes its implementation on a real robot challenging \cite{mystkowski2013robust}.

\section{Conclusion}
This research has studied the robust attitude control of a quadrotor with taking into account the uncertainty and inputs constraints. Knowing that the model and parameters of the system do not exist, a continuous-time black-box model identification based on the real sampled data is adopted to acquire both nominated linear model of the system and the uncertainty model encapsulating the deviation of the nonlinear system from the nominal model. Now, the essential prerequisites for applying the robust $H_{\infty}$ controller is provided. By solving the mixed-sensitivity problem, a robust stabilizing $H_{\infty}$ controller is obtained which satisfies tracking, disturbance attenuation and input saturation objectives.
Having accomplished that, the pitch and roll controllers are implemented on the robot, and qualitative and quantitative experiments are performed to verify the performance of $H_{\infty}$ controllers. According to both experimental and simulation results, the controllers successfully fulfill the aforementioned objectives. Finally, by comparing the performance of the proposed controllers with PID and $\mu$ synthesis controllers, it is demonstrated that the $H_{\infty}$ controller has a better performance.


\bibliography{RefAbbv}

\begin{thebibliography}{10}
\expandafter\ifx\csname url\endcsname\relax
  \def\url#1{\texttt{#1}}\fi
\expandafter\ifx\csname urlprefix\endcsname\relax\def\urlprefix{URL }\fi
\expandafter\ifx\csname href\endcsname\relax
  \def\href#1#2{#2} \def\path#1{#1}\fi

\bibitem{wang2016detecting}
L.~Wang, F.~Chen, H.~Yin, Detecting and tracking vehicles in traffic by
  unmanned aerial vehicles, Automat. Constr. 72 (2016) 294--308.

\bibitem{berra2017commercial}
E.~F. Berra, R.~Gaulton, S.~Barr, Commercial off-the-shelf digital cameras on
  unmanned aerial vehicles for multitemporal monitoring of vegetation
  reflectance and ndvi, IEEE Trans. Geosci. Remote Sens. 55~(9) (2017)
  4878--4886.

\bibitem{acevedo2014decentralized}
J.~J. Acevedo, B.~C. Arrue, I.~Maza, A.~Ollero, A decentralized algorithm for
  area surveillance missions using a team of aerial robots with different
  sensing capabilities, in: Robot. Autom. (ICRA), 2014 IEEE Int. Conf., IEEE,
  2014, pp. 4735--4740.

\bibitem{kingston2008decentralized}
D.~Kingston, R.~W. Beard, R.~S. Holt, Decentralized perimeter surveillance
  using a team of {UAV}s, IEEE Trans. Robot. 24~(6) (2008) 1394--1404.

\bibitem{silvagni2017multipurpose}
M.~Silvagni, A.~Tonoli, E.~Zenerino, M.~Chiaberge, Multipurpose {UAV} for
  search and rescue operations in mountain avalanche events, Geomat. Nat. Haz.
  Risk. 8~(1) (2017) 18--33.

\bibitem{khan2015information}
A.~Khan, E.~Yanmaz, B.~Rinner, Information exchange and decision making in
  micro aerial vehicle networks for cooperative search, IEEE Trans. Control
  Netw. Syst. 2~(4) (2015) 335--347.

\bibitem{tokekar2016sensor}
P.~Tokekar, J.~Vander~Hook, D.~Mulla, V.~Isler, Sensor planning for a symbiotic
  {UAV} and {UGV} system for precision agriculture, IEEE Trans. Robot 32~(6)
  (2016) 1498--1511.

\bibitem{tetila2017identification}
E.~C. Tetila, B.~B. Machado, N.~A. de~Souza~Belete, D.~A. Guimar{\~a}es,
  H.~Pistori, Identification of soybean foliar diseases using unmanned aerial
  vehicle images, IEEE Geosci. Remote Sens. Lett. 14~(12) (2017) 2190--2194.

\bibitem{gawel2017aerial}
A.~Gawel, M.~Kamel, T.~Novkovic, J.~Widauer, D.~Schindler, B.~P. von
  Altishofen, R.~Siegwart, J.~Nieto, Aerial picking and delivery of magnetic
  objects with mavs, in: Robot. Autom. (ICRA), 2017 IEEE Int. Conf., IEEE,
  2017, pp. 5746--5752.

\bibitem{arbanas2016aerial}
B.~Arbanas, A.~Ivanovic, M.~Car, T.~Haus, M.~Orsag, T.~Petrovic, S.~Bogdan,
  Aerial-ground robotic system for autonomous delivery tasks, in: Robot. Autom.
  (ICRA), 2016 IEEE Int. Conf., IEEE, 2016, pp. 5463--5468.

\bibitem{hoffmann2007quadrotor}
G.~Hoffmann, H.~Huang, S.~Waslander, C.~Tomlin, Quadrotor helicopter flight
  dynamics and control: Theory and experiment, in: AIAA Guid., Nav. Control
  Conf. Exhibit, 2007, p. 6461.

\bibitem{szafranski2011different}
G.~Szafranski, R.~Czyba, Different approaches of {PID} control {UAV} type
  quadrotor, in: Proc. Int. Micro Air Veh. conf., 2011.

\bibitem{garcia2012robust}
R.~Garcia, F.~Rubio, M.~Ortega, Robust {PID} control of the quadrotor
  helicopter, IFAC Proc. Vols. 45~(3) (2012) 229--234.

\bibitem{zhou2017robust}
J.~Zhou, R.~Deng, Z.~Shi, Y.~Zhong, Robust cascade {PID} attitude control of
  quadrotor helicopters subject to wind disturbance, in: Control Conf. (CCC),
  2017 36th Chinese, IEEE, 2017, pp. 6558--6563.

\bibitem{bouabdallah2004pid}
S.~Bouabdallah, A.~Noth, R.~Siegwart, {PID} vs {LQ} control techniques applied
  to an indoor micro quadrotor, in: Proc. IEEE Int. Conf. Intell. Robots Syst.
  (IROS), IEEE, 2004, pp. 2451--2456.

\bibitem{liu2013robust}
H.~Liu, D.~Derawi, J.~Kim, Y.~Zhong, Robust optimal attitude control of
  hexarotor robotic vehicles, Nonlinear Dynam. 74~(4) (2013) 1155--1168.

\bibitem{flores2013lyapunov}
G.~Flores, R.~Lozano, Lyapunov-based controller using singular perturbation
  theory: An application on a mini-uav, in: Amer. Control Conf. (ACC), 2013,
  IEEE, 2013, pp. 1596--1601.

\bibitem{huo2014attitude}
X.~Huo, M.~Huo, H.~R. Karimi, Attitude stabilization control of a quadrotor
  {UAV} by using backstepping approach, Math. Probl. Eng. 2014.

\bibitem{das2009backstepping}
A.~Das, F.~Lewis, K.~Subbarao, Backstepping approach for controlling a
  quadrotor using lagrange form dynamics, J. Intell. Robot. Syst. 56~(1-2)
  (2009) 127--151.

\bibitem{djamel2016attitude}
K.~Djamel, M.~Abdellah, A.~Benallegue, Attitude optimal backstepping controller
  based quaternion for a {UAV}, Math. Probl. Eng. 2016.

\bibitem{shao2018robust}
X.~Shao, J.~Liu, H.~Wang, Robust back-stepping output feedback trajectory
  tracking for quadrotors via extended state observer and sigmoid tracking
  differentiator, Mech. Syst. Signal Process. 104 (2018) 631--647.

\bibitem{chen2016robust}
F.~Chen, R.~Jiang, K.~Zhang, B.~Jiang, G.~Tao, Robust backstepping sliding-mode
  control and observer-based fault estimation for a quadrotor {UAV}, IEEE
  Trans. Ind. Electron. 63~(8) (2016) 5044--5056.

\bibitem{jia2017integral}
Z.~Jia, J.~Yu, Y.~Mei, Y.~Chen, Y.~Shen, X.~Ai, Integral backstepping sliding
  mode control for quadrotor helicopter under external uncertain disturbances,
  Aerosp. Sci. Technol. 68 (2017) 299--307.

\bibitem{yang2016attitude}
Y.~Yang, Y.~Yan, Attitude regulation for unmanned quadrotors using adaptive
  fuzzy gain-scheduling sliding mode control, Aerosp. Sci. Technol. 54 (2016)
  208--217.

\bibitem{zheng2014second}
E.-H. Zheng, J.-J. Xiong, J.-L. Luo, Second order sliding mode control for a
  quadrotor {UAV}, ISA Trans. 53~(4) (2014) 1350--1356.

\bibitem{wang2019disturbance}
B.~Wang, X.~Yu, L.~Mu, Y.~Zhang, Disturbance observer-based adaptive
  fault-tolerant control for a quadrotor helicopter subject to parametric
  uncertainties and external disturbances, Mech. Syst. Signal Process. 120
  (2019) 727--743.

\bibitem{mallavalli2018fault}
S.~Mallavalli, A.~Fekih, A fault tolerant tracking control for a quadrotor
  {UAV} subject to simultaneous actuator faults and exogenous disturbances,
  Int. J. Control (2018) 1--14.

\bibitem{liu2017robust}
H.~Liu, W.~Zhao, Z.~Zuo, Y.~Zhong, Robust control for quadrotors with multiple
  time-varying uncertainties and delays, IEEE Trans. Ind. Electron. 64~(2)
  (2017) 1303--1312.

\bibitem{raffo2010integral}
G.~V. Raffo, M.~G. Ortega, F.~R. Rubio, An integral predictive/nonlinear
  {H}$\infty$ control structure for a quadrotor helicopter, Automatica 46~(1)
  (2010) 29--39.

\bibitem{kerma2012nonlinear}
M.~Kerma, A.~Mokhtari, B.~Abdelaziz, Y.~Orlov, Nonlinear {H}$\infty$ control of
  a quadrotor ({UAV}), using high order sliding mode disturbance estimator,
  Int. J. Control 85~(12) (2012) 1876--1885.

\bibitem{hu2018adaptive}
Y.~Hu, W.~Gu, H.~Zhang, H.~Chen, Adaptive robust triple-step control for
  compensating cogging torque and model uncertainty in a dc motor, IEEE Trans.
  Syst., Man, Cybern., Syst.

\bibitem{zhang2016active}
H.~Zhang, J.~Wang, Active steering actuator fault detection for an
  automatically-steered electric ground vehicle, IEEE Trans. Veh. Technol.
  66~(5) (2016) 3685--3702.

\bibitem{jiang2019hydrothermal}
K.~Jiang, F.~Yan, H.~Zhang, Hydrothermal aging factor estimation for two-cell
  diesel-engine scr systems via a dual time-scale unscented kalman filter, IEEE
  Trans. Ind. Electron.

\bibitem{tran2018adaptive}
T.-T. Tran, S.~S. Ge, W.~He, Adaptive control of a quadrotor aerial vehicle
  with input constraints and uncertain parameters, Int. J. Control 91~(5)
  (2018) 1140--1160.

\bibitem{liu2018parameter}
W.~Liu, X.~Huo, J.~Liu, L.~Wang, Parameter identification for a quadrotor
  helicopter using multivariable extremum seeking algorithm, Int. J. Control.
  Autom. Syst. 16~(4) (2018) 1951--1961.

\bibitem{bergamasco2014identification}
M.~Bergamasco, M.~Lovera, Identification of linear models for the dynamics of a
  hovering quadrotor, IEEE Trans. Control Syst. Technol. 22~(5) (2014)
  1696--1707.

\bibitem{agand2017decentralized}
P.~Agand, M.~Motaharifar, H.~D. Taghirad, Decentralized robust control for
  teleoperated needle insertion with uncertainty and communication delay,
  Mechatronics 46 (2017) 46--59.

\bibitem{bataleblu2016robust}
A.~Bataleblu, M.~Motaharifar, E.~Abedlu, H.~D. Taghirad, Robust {H}$\infty$
  control of a 2rt parallel robot for eye surgery, in: Robot. Mechatronics
  (ICROM), 2016 4th Int. Conf., IEEE, 2016, pp. 136--141.

\bibitem{taghirad2001h}
H.~Taghirad, P.~Belanger, {H}$\infty$-based robust torque control of harmonic
  drive systems, J. Dyn. Sys. Meas. Control 123~(3) (2001) 338--345.

\bibitem{golovin2019robust}
I.~Golovin, S.~Palis, Robust control for active damping of elastic gantry crane
  vibrations, Mech. Syst. Signal Process. 121 (2019) 264--278.

\bibitem{mystkowski2016mu}
A.~Mystkowski, A.~P. Koszewnik, Mu-synthesis robust control of 3{D} bar
  structure vibration using piezo-stack actuators, Mech. Syst. Signal Process.
  78 (2016) 18--27.

\bibitem{safaee2013system}
A.~Safaee, H.~D. Taghirad, System identification and robust controller design
  for the autopilot of an unmanned helicopter, in: Control Conf. (ASCC), 2013
  9th Asian, IEEE, 2013, pp. 1--6.

\bibitem{wang2013robust}
X.~Wang, G.~Lu, Y.~Zhong, Robust {H}$\infty$ attitude control of a laboratory
  helicopter, Robot. Auton. Syst. 61~(12) (2013) 1247--1257.

\bibitem{mahony2008nonlinear}
R.~Mahony, T.~Hamel, J.~M. Pflimlin, Nonlinear complementary filters on the
  special orthogonal group, IEEE Trans. Autom. Control 53~(5) (2008)
  1203--1218.

\bibitem{mystkowski2013robust}
A.~Mystkowski, Robust control of the micro {UAV} dynamics with an autopilot, J.
  Theor. Appl. Mech. 51~(3) (2013) 751--761.

\end{thebibliography}

\end{document}